\documentclass[a4paper]{article}

\usepackage{amsthm,amsmath,amssymb}
\usepackage{enumitem,listings}
\usepackage{algorithm,algpseudocode}
\usepackage{multirow}

\usepackage{tikz,graphicx}
\usepackage{standalone}

\usepackage[title]{appendix}

\usepackage[hidelinks]{hyperref}
\usepackage{cleveref}

\theoremstyle{theorem}
\newtheorem{theorem}{Theorem}[section]

\newtheorem{corollary}[theorem]{Corollary}
\newtheorem{lemma}[theorem]{Lemma}
\newtheorem{conjecture}[theorem]{Conjecture}

\newtheorem{algorithm_thm}[theorem]{Algorithm}

\theoremstyle{definition}
\newtheorem{definition}[theorem]{Definition}

\counterwithin{figure}{section}
\counterwithin{table}{section}

\crefname{conjecture}{conjecture}{conjectures}

\usetikzlibrary{positioning}
\usetikzlibrary{decorations.markings}
\usetikzlibrary{arrows}

\newcommand{\tri}{\mathcal{T}}
\newcommand{\manifold}{\mathcal{M}}
\newcommand{\CPplusfamily}[1]{\mathcal{P}(#1)} 
\newcommand{\CPplusminusfamily}[1]{\mathcal{A}(#1)} 
\newcommand{\SSfamily}[1]{\mathcal{E}(#1)} 
\newcommand{\CPSSfamily}[2]{\mathcal{D}(#1,#2)} 

\newcommand{\id}{\textrm{id}}

\newcommand{\NN}{\mathbb{N}}
\newcommand{\ZZ}{\mathbb{Z}}

\newcommand{\CP}{\mathbb{CP}}
\newcommand{\Sym}{\operatorname{Sym}}

\DeclareMathOperator{\rank}{rank}
\DeclareMathOperator{\sign}{sign}
\DeclareMathOperator{\openstar}{star}
\DeclareMathOperator{\link}{link}

\begin{document}

\title{Small Triangulations of Simply Connected $4$-Manifolds}

\author{Jonathan Spreer and Lucy Tobin}
\date{}

\maketitle

\begin{abstract}
  We present small triangulations of all connected sums of $\CP^2$ and $S^2 \times S^2$ with the standard piecewise linear structure. Our triangulations have $2\beta_2+2$ pentachora, where $\beta_2$ is the second Betti number of the manifold. By a conjecture of the authors and, independently, Burke, these triangulations have the smallest possible number of pentachora for their respective topological types.
\end{abstract}



\section{Introduction}

Given a piecewise linear manifold $\manifold$, an important measure for how ``complicated'' it is topologically, is the smallest number of simplices needed to build $\manifold$ from gluing these simplices along their faces, that is, from {\em triangulating} $\manifold$. In the literature, this quantity is referred to as the {\em complexity} $c(\manifold)$ of $\manifold$ \cite{matveev90-complexity}. 

A triangulation realising the complexity of its underlying manifold is said to be {\em minimal}. Naturally, the size of a minimal triangulation depends on the exact rules of how we are allowed to glue simplices together. If the triangulation is required to be a simplicial complex, then the complexity is larger. If graph-encodings of manifolds are considered (that is, regular cell complexes, that allow a rainbow colouring of their vertices), a minimal triangulation is typically considerably smaller. Finally, the truly smallest number of simplices needed to triangulate a manifold is achieved for what is called {\em generalised triangulations}, where even pairs of faces from the same simplex can be identified, and a triangulation can have as little as one vertex.

\begin{figure}
\centering
\resizebox{0.9\linewidth}{!}{\tikzset{
    >=latex,
    singlearrow/.style={decoration={markings, mark=at position 0.5*\pgfdecoratedpathlength+4pt with {\arrow{#1}}}, postaction={decorate}},
    doublearrow/.style={decoration={markings, mark=at position 0.5*\pgfdecoratedpathlength+7pt with {\arrow{{#1}{#1}}}}, postaction={decorate}},
    triplearrow/.style={decoration={markings, mark=at position 0.5*\pgfdecoratedpathlength+9pt with {\arrow{{#1}{#1}{#1}}}}, postaction={decorate}}
}

\begin{tabular}{c@{\hspace{2cm}}c@{\hspace{2cm}}c}

\begin{tikzpicture}[every node/.style={circle, inner sep=0, outer sep=0, minimum width=0.1cm, fill=black}, every path/.style={thick}, label distance=2mm, on grid]

\node (11) at (0,0) {} ;
\node (12) [right=of 11] {} ;
\node (13) [right=of 12] {} ;
\node (14) [right=of 13] {} ;

\node (21) [above=of 11] {} ;
\node (22) [right=of 21] {} ;
\node (24) [above=of 14] {} ;

\node (31) [above=of 21] {} ;
\node (33) [right=2 of 31] {} ;
\node (34) [above=of 24] {} ;

\node (41) [above=of 31] {} ;
\node (42) [right=of 41] {} ;
\node (43) [right=of 42] {} ;
\node (44) [right=of 43] {} ;

\draw[singlearrow=>] (11) -- (12) ;
\draw[doublearrow=>] (12) -- (13) ;
\draw[triplearrow=>] (13) -- (14) ;

\draw[singlearrow=>] (41) -- (42) ;
\draw[doublearrow=>] (42) -- (43) ;
\draw[triplearrow=>] (43) -- (44) ;

\draw[decoration={markings, mark=at position 0.5*\pgfdecoratedpathlength+1pt with {\arrow{To}}}, postaction={decorate}] (11) -- (21) ;
\draw[decoration={markings, mark=at position 0.5*\pgfdecoratedpathlength+3.5pt with {\arrow{{To}{To}}}}, postaction={decorate}] (21) -- (31) ;
\draw[decoration={markings, mark=at position 0.5*\pgfdecoratedpathlength+6pt with {\arrow{{To}{To}{To}}}}, postaction={decorate}] (31) -- (41) ;

\draw[decoration={markings, mark=at position 0.5*\pgfdecoratedpathlength+1pt with {\arrow{To}}}, postaction={decorate}] (14) -- (24) ;
\draw[decoration={markings, mark=at position 0.5*\pgfdecoratedpathlength+3.5pt with {\arrow{{To}{To}}}}, postaction={decorate}] (24) -- (34) ;
\draw[decoration={markings, mark=at position 0.5*\pgfdecoratedpathlength+6pt with {\arrow{{To}{To}{To}}}}, postaction={decorate}] (34) -- (44) ;

\draw (12) -- (21) ;
\draw (13) -- (24) ;
\draw (42) -- (31) ;
\draw (43) -- (34) ;

\draw (12) -- (24) ;
\draw (31) -- (43) ;

\draw (21) -- (22) ;
\draw (31) -- (22) ;
\draw (24) -- (33) ;
\draw (34) -- (33) ;

\draw (12) -- (22) ;
\draw (22) -- (24) ;
\draw (22) -- (33) ;
\draw (31) -- (33) ;
\draw (33) -- (43) ;

\end{tikzpicture}

&

\begin{tikzpicture}[every node/.style={circle, inner sep=0, outer sep=0, minimum width=0.1cm, fill=black}, every path/.style={thick}, label distance=2mm, on grid]

\node[red] (c) at (0,0) {} ;
\path +(-120:{3/sqrt(3)}) node[blue] (11) {} ;
\path +(-60:{3/sqrt(3)}) node[green] (12) {} ;
\path +(180:{3/sqrt(3)}) node[green] (21) {} ;
\path +(0:{3/sqrt(3)}) node[blue] (22) {} ;
\path +(120:{3/sqrt(3)}) node[blue] (31) {} ;
\path +(60:{3/sqrt(3)}) node[green] (32) {} ;

\node[fill=none] (spacertop) at (0,1.5) {} ;
\node[fill=none] (spacerbot) at (0,-1.5) {} ;

\draw (c) -- (11) ;
\draw (c) -- (12) ;
\draw (c) -- (21) ;
\draw (c) -- (22) ;
\draw (c) -- (31) ;
\draw (c) -- (32) ;

\draw[singlearrow=>] (11) -- (12) ;
\draw[doublearrow={>}] (12) -- (22) ;
\draw[triplearrow=>] (22) -- (32) ;
\draw[singlearrow=>] (31) -- (32) ;
\draw[doublearrow=>] (21) -- (31) ;
\draw[triplearrow=>] (11) -- (21) ;

\end{tikzpicture}

&

\begin{tikzpicture}[every node/.style={circle, inner sep=0, outer sep=0, minimum width=0.1cm, fill=black}, every path/.style={thick}, label distance=2mm, on grid]

\node (11) at (0,0) {} ;
\node (12) [right=3 of 11] {} ;

\node (21) [above=3 of 11] {} ;
\node (22) [right=3 of 21] {} ;

\draw[singlearrow=>] (11) -- (12) ;
\draw[singlearrow=>] (21) -- (22) ;
\draw[doublearrow=>] (11) -- (21) ;
\draw[doublearrow=>] (12) -- (22) ;

\draw (11) -- (22) ;

\end{tikzpicture}

\\

\\

\begin{tikzpicture}[every node/.style={circle, inner sep=0, outer sep=0, minimum width=0.1cm, fill=black}, every path/.style={thick}, label distance=2mm, on grid]

\node (11) at (0,0) {} ;
\node (12) [right=3/4 of 11] {} ;
\node (13) [right=3/4 of 12] {} ;
\node (14) [right=3/4 of 13] {} ;
\node (15) [right=3/4 of 14] {} ;

\node (21) at (3/8,1) {} ;
\node (22) [right=3/4 of 21] {} ;
\node (23) [right=3/4 of 22] {} ;
\node (24) [right=3/4 of 23] {} ;

\node (31) [above=2 of 11] {} ;
\node (32) [right=3/4 of 31] {} ;
\node (33) [right=3/4 of 32] {} ;
\node (34) [right=3/4 of 33] {} ;
\node (35) [right=3/4 of 34] {} ;

\draw (11) -- (12) ;
\draw (12) -- (13) ;
\draw (13) -- (14) ;
\draw (14) -- (15) ;

\draw (21) -- (22) ;
\draw (22) -- (23) ;
\draw (23) -- (24) ;

\draw (31) -- (32) ;
\draw (32) -- (33) ;
\draw (33) -- (34) ;
\draw (34) -- (35) ;

\draw (12) -- (21) ;
\draw (13) -- (23) ;

\draw (22) -- (33) ;
\draw (24) -- (34) ;

\draw (11) to[bend right=30] (15) ;
\draw (21) to[bend left=30] (24) ;
\draw (31) to[bend left=30] (35) ;

\draw (11) -- (31) ;
\draw (14) -- (32) ;
\draw (15) -- (35) ;

\end{tikzpicture}

&

\begin{tikzpicture}[every node/.style={circle, inner sep=0, outer sep=0, minimum width=0.1cm, fill=black}, every path/.style={thick}, label distance=2mm, on grid]

\node[fill=none] (c) at (0,0) {} ;
\path +(-90:1.5) node (11) {} ;
\path +(-150:1.5) node (21) {} ;
\path +(-30:1.5) node (22) {} ;
\path +(150:1.5) node (31) {} ;
\path +(30:1.5) node (32) {} ;
\path +(90:1.5) node (41) {} ;

\draw[blue] (11) -- (22) ;
\draw[green] (22) -- (32) ;
\draw[blue] (32) -- (41) ;
\draw[green] (41) -- (31) ;
\draw[blue] (31) -- (21) ;
\draw[green] (21) -- (11) ;

\draw[red] (11) -- (41) ;
\draw[red] (21) -- (32) ;
\draw[red] (31) -- (22) ;

\end{tikzpicture}

&

\begin{tikzpicture}[every node/.style={circle, inner sep=0, outer sep=0, minimum width=0.1cm, fill=black}, every path/.style={thick}, label distance=2mm, on grid]

\node[fill=none] (c) at (0,0) {} ;
\node (tl) [above left=of c] {} ;
\node (br) [below right=of c] {} ;

\node[fill=none] (spacertop) at (0,1.5) {} ;
\node[fill=none] (spacerbot) at (0,-1.5) {} ;

\draw (tl) -- (br) ;
\draw (tl) to[bend left=60] (br) ;
\draw (tl) to[bend right=60] (br) ;

\end{tikzpicture}

\end{tabular}}
\caption{Minimal triangulations of the torus as a simplicial complex (left), graph encoding (middle), and generalised triangulation (right). Top row shows triangulations and bottom row shows corresponding dual graphs.\label{fig:torus}}
\end{figure}
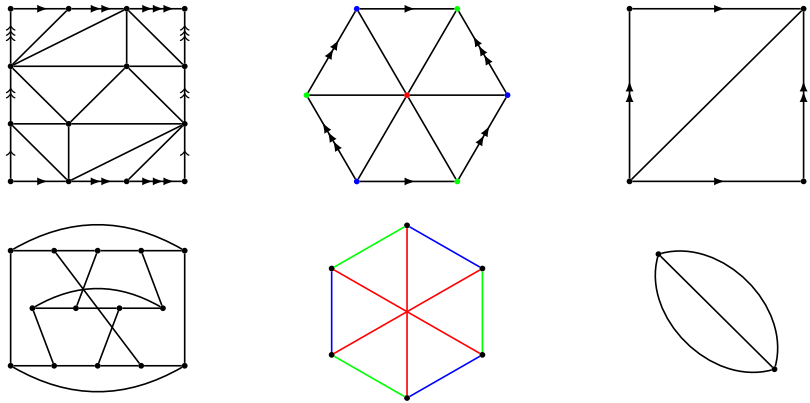

For instance, the smallest simplicial complex triangulating the torus has 14 triangles, the smallest graph-encoding has six, and the smallest generalised triangulation has only two, see \Cref{fig:torus}.

Minimal triangulations are known for all surfaces in all three settings, see, for instance, \cite{jaco06-layered,Jungermann80Surfaces,Shankar23Sampling}. For generalised triangulations, minimal triangulations of a surface $S$ of non-positive Euler characteristic $\chi(S)$ are precisely those with a single vertex, and $S$ has complexity $2 - 2\chi(S)$. This holds for both orientable and non-orientable surfaces.

For triangulations of $3$-dimensional manifolds, the situation is very different. Bounds for the complexity are known for some simplicial triangulations of $3$-manifolds \cite{lutz09-fvectors,sulanke09-enumeration}. For graph-encodings, infinite families of minimal triangulations exist \cite{Basak14MinTrig}. For generalised triangulations, several techniques are known to obtain bounds and exact complexity results for infinite families of manifolds, see \cite{Cha-complexities-2018, JJSTBoundsOnNormSurfs,jaco09-minimal-lens,Jaco-ideal-2020,DehnFilling,Lackenby-2019-complexity,Rubinstein2024} for a collection of results. However, most of these results only cover particular classes of manifolds, and practical bounds remain out of reach in the general case.

In dimension four, results in the literature are largely focused on simply connected $4$-manifolds. There are minimal simplicial triangulations of the complex projective plane $\mathbb{C}P^2$ \cite{Banchoff83CP2}, and a manifold homeomorphic to the K3 surface \cite{Casella01K3}. In terms of graph encodings, in \cite{BasakSpreer2016} Basak together with the first author presented smallest possible triangulations with the minimal number of vertices and edges for all standard simply connected $4$-manifolds (modulo the $11/8$-conjecture, see \Cref{conj:11/8}). 

Finally, in arbitrary higher dimensions, some results exist in the simplicial setting. Specifically, see \cite{KuehnelBook} for minimal triangulations of circle bundles over the sphere, \cite{KleeCentrallySymmetric} for minimal triangulations of sphere bundles in the centrally symmetric setting, and \cite{Kuehnel99CensusTight} for a conjecture relating minimal to what are called {\em tight} triangulations.

\medskip

In this article, we present generalised triangulations of simply connected $4$-manifolds with second Betti number $k$, built from $2k+2$ four-dimensional simplices. Our main result is summarised in the following theorem.

\begin{theorem}
  \label{thm:introduction}
  Let $\manifold$ be a simply connected $4$-manifold obtained from connected sums of PL standard summands $\mathbb{C}P^2$ and $S^2 \times S^2$ in any orientation, then $\manifold$ can be triangulated using  $2 \beta_2(\manifold) + 2$ pentachora.
\end{theorem}

To the best of our knowledge, these triangulations are the smallest known triangulations of the simply connected manifolds in question. In fact, we have the following conjecture, first stated by the authors in \cite{SpreerTobin24FaceNumbers} and independently by Burke in \cite{Burke2024}. 

\begin{conjecture}
  Let $\manifold$ be a simply connected PL $4$-manifold, then the complexity of $\manifold$ satisfies $c(\manifold) \geq 2 \beta_2 (\manifold) + 2$.
\end{conjecture}

Note that we do not have similar triangulations of simply connected $4$-manifolds containing the K3 surface as connected summands. In fact, no triangulation of the K3 surface with $2 \beta_2 + 2 = 46$ four-simplices is know to date. The smallest known triangulation of the K3 surface has $54$ simplices and is due to Burke \cite{Burke2024}. 

\medskip

This paper is organised as follows: after going through some necessary background in \Cref{sec:background}, we present the setup for our infinite families of triangulations in \Cref{sec:setup}. In \Cref{sec:series} we present our actual triangulations before identifying their PL-homeomorphism type in \Cref{sec:PLtype}. Finally, we describe the implementation of our computer-aided proof methods in \Cref{sec:computation}. We list sequences verifying the PL-homeomorphism type of our triangulations in \Cref{app:sequences}, and code to generate our triangulations in \Cref{app:code}.

\paragraph*{Acknowledgements.}

The work of Spreer is partially supported by the Australian Research Council under the Discovery Project scheme (grant number DP220102588). Part of this work was done while both authors were visiting Technische Universit\"at Berlin. Our special thanks go to Michael Joswig and his research group for their support and hospitality. We also want to thank the Discrete Geometry group at Freie Universt\"at Berlin for additional support.

\section{Background}
\label{sec:background}

\subsection{4-Manifolds}
\label{subsec:manifold}

This section covers basic results on 4-manifolds and their intersection forms, and may be skipped by readers familiar with 4-manifold topology. For a more thorough introduction into 4-manifold topology, we refer to \cite{GompfStipsicz1999}.

When working in dimension 4, one must distinguish between \textit{topological} manifolds, and \textit{smooth} and \textit{piecewise-linear (PL)} manifolds: While for every smooth structure on a 4-manifold, there exists a corresponding PL structure, and vice versa, a single compact topological 4-manifold may admit between $0$ and countably many distinct smooth (or PL) structures.

In this article, we focus on smooth 4-manifolds and \textit{smoothable} topological 4-manifolds -- that is, topological 4-manifolds that admit at least one smooth structure. A smoothable topological 4-manifold can be thought of as a (non-empty) collection of smooth 4-manifolds that are pairwise homeomorphic, but not diffeomorphic.

The \textit{intersection form} $Q_\manifold: H^2(\manifold) \times H^2(\manifold) \to \ZZ$ of a simply connected, oriented, closed 4-manifold $\manifold$ is a unimodular, symmetric, bilinear integral form which maps a pair of $2$-cohomology classes to their {\em cup-product}, roughly describing their number of intersections, see \cite{GompfStipsicz1999} for a precise definition. The intersection form has a \textit{rank} $\rank(Q_\manifold)=\rank(H^2(\manifold)) \in \NN$, a \textit{signature} $\sign(Q_\manifold) \in \ZZ$, and a \textit{type}, even or odd. If $\sign(Q_\manifold)=\rank(Q_\manifold)$ ($\sign(Q_\manifold)=-\rank(Q_\manifold)$) it is called \textit{positive definite} (resp. \textit{negative definite}), and \textit{indefinite} otherwise.

Although the intersection form is defined on \textit{oriented} manifolds, we will mostly consider \textit{orientable} manifolds (every simply connected 4-manifold is necessarily orientable) with the observation that reversing a chosen orientation simply negates the form: $Q_{\overline{\manifold}} = -Q_\manifold$ and $\sign(Q_{\overline{\manifold}})=-\sign(Q_\manifold)$, where $\overline{\manifold}$ denotes $\manifold$ with opposite orientation.

We now briefly recall fundamental results, required for \Cref{subsec:cp2-s2xs2-csums}, beginning with Freedman's famous classification theorem:

\begin{theorem}[Freedman \cite{Freedman1982}]
  Every unimodular, symmetric, bilinear integral form is the intersection form of a simply connected, closed, topological 4-manifold. Moreover, (up to homeomorphism of manifolds and isomorphism of forms):
  \begin{itemize}
    \item if the form is even, there is exactly one such 4-manifold
    \item if the form is odd, there are exactly two such manifolds and at least one does not admit a smooth structure
  \end{itemize}
  \label{thm:freedman}
\end{theorem}

In particular, the topological type of a simply connected, closed, \textit{smooth} 4-manifold is uniquely determined by its intersection form. Which forms are the intersection form of a smoothable 4-manifold is restricted by the following two results (and one conjecture):

\begin{theorem}[Donaldson \cite{Donaldson1983}]
  Let $\manifold$ be a simply connected, closed (oriented) smooth 4-manifold with definite intersection form $Q$. Then $Q \cong k \cdot (+1)$ (if positive definite) or $Q \cong k \cdot (-1)$ (if negative definite) for some $k \in \NN$.
  \label{thm:donaldson}
\end{theorem}

\begin{theorem}[Rokhlin \cite{Rokhlin1952}]
  Let $\manifold$ be a simply connected, closed, smooth 4-manifold with even intersection form $Q$. Then $\sign(Q) \equiv 0 \textrm{ mod } 16$.
  \label{thm:rokhlin}
\end{theorem}

\begin{conjecture}[The 11/8 Conjecture]
  Let $\manifold$ be a simply connected, closed, smooth 4-manifold with even intersection form $Q$. Then $\rank(Q) \ge 11/8 \cdot |\sign(Q)|$.
  \label{conj:11/8}
\end{conjecture}

\hyperref[thm:donaldson]{Donaldson's theorem} covers the case of definite intersection forms. On the other hand, we have the following results for indefinite forms:

\begin{lemma}[\cite{Milnor2013}]
  Let $Q$ be an odd and indefinite unimodular, symmetric, bilinear form. Then
  \begin{align*}
    Q \cong \frac{\rank(Q)+\sign(Q)}{2} \cdot (+1) \oplus \frac{\rank(Q)-\sign(Q)}{2} \cdot (-1)
  \end{align*}
  \label{lem:odd-indefinite-forms}
\end{lemma}

\begin{theorem}[\cite{Milnor2013}]
  An indefinite unimodular, symmetric, bilinear form is determined up to isomorphism by its rank, signature and type.
  \label{thm:indefinite-forms-classification}
\end{theorem}

\hyperref[thm:freedman]{Freedman's theorem}, \hyperref[thm:donaldson]{Donaldson's theorem} and \Cref{thm:indefinite-forms-classification} together imply:

\begin{corollary}
  The homeomorphism type of a simply connected, closed, smooth 4-manifold is uniquely determined by the rank, signature and type of its intersection form.
  \label{cor:rank-sig-type}
\end{corollary}

\subsection{Triangulations}
\label{sec:triangulation}

In this work we consider \textit{generalised triangulations} which, in contrast to simplicial complexes, allow for two simplices to be glued together across multiple different faces, and for two faces of the same simplex to be glued together.

A \textit{(generalised) triangulation of a $d$-manifold} $\tri$ is a set of abstract $d$-simplices, called \textit{facets}, together with a set of affine bijections which identify two distinct $(d-1)$-dimensional faces of the facets, called \textit{gluings}. As a result of these gluings, faces become identified and we refer to a resulting class of faces as a single face of $\tri$. The number of pre-glued faces in this class is referred to as the {\em degree} of the face. No face can be identified with itself along a non-identity mapping. Faces of dimensions $0$, $1$, $2$, $3$, $4$, and $(d-1)$ are referred to as {\em vertices}, {\em edges}, {\em triangles}, {\em tetrahedra}, {\em pentachora}, and {\em ridges} of $\tri$, respectively. The boundary of a small neighbourhood of a vertex of $\tri$, its {\em link}, see \Cref{def:starlink}, must triangulate a $(d-1)$-dimensional sphere with standard PL structure. The definition is completed by noting that a $0$-dimensional sphere consists of two isolated vertices. We say that $\tri$ {\em triangulates} a PL $d$-manifold $\manifold$, if its underlying set $|\tri|$ is PL-homeomorphic to $\manifold$. 

Note that our definition of a triangulation implies that $\manifold$ is closed, that is, every ridge appears in exactly one gluing. By admitting vertex links to be triangulated $(d-1)$-balls as well, we can extend this definition to triangulations of manifolds with non-empty boundary. In this case the {\em boundary} of the triangulation $\tri$, denoted by $\partial \tri$, is comprised of all un-glued ridges together with their faces. In this case, all faces that are not in the boundary are referred to as {\em internal faces}. Furthermore, we will sometimes need to consider triangulations with no conditions on the vertex links, where $|\tri|$ is not necessarily a manifold. We refer to these simply as {\em triangulations}.

The total numbers of $i$-dimensional faces (or simply \textit{$i$-faces}) $f_i$ of a $d$-dimensional triangulation $\tri$, $0 \leq i \leq d$, together form the \textit{$f$-vector} of $\tri$, written $f(\tri)=(f_0,f_1,\dots,f_d)$. 

Note that, before gluing, a $(d-1)$-face of a facet is determined by a set of $d$ labelled vertices of that facet, and a gluing is completely determined by a bijection between two such sets. As such, we consider the vertices of a facet to be labelled by $\{0,1,\dots,d\}$, and encode a gluing as a tuple $(\Delta_1,\Delta_2,i,\sigma)$, where $\Delta_1$, $\Delta_2$ are facets, $i \in \{0,1,\dots,d\}$ is the vertex of $\Delta_1$ \textit{not} involved in the gluing, $\sigma \in \Sym_{d+1}$ is a permutation such that $\sigma(i)$ is the vertex of $\Delta_2$ not involved in the gluing, and vertex $j$ of $\Delta_1$ is matched with vertex $\sigma(j)$ of $\Delta_2$ for each $j \in \{0,1,\dots,d\} \setminus \{i\}$.

Two triangulations are \textit{combinatorially equivalent}, or {\em isomorphic}, if one can be obtained from the other by relabelling its facets and vertices. We are typically interested in triangulations up to combinatorial equivalence, but must work with specific labelled triangulations to describe operations on them.

\begin{definition}
  \label{def:starlink}
  Let $F$ be an $i$-face in a $d$-dimensional triangulation $\tri$. The \textit{(open) star} $\openstar(F)$ is the set of all faces of $\tri$ which have $F$ as a face. The \textit{link} $\link(F)$ is the $(d-i-1)$-dimensional triangulation constructed as follows:
  \begin{itemize}
    \item For each facet $\Delta$ containing $F$, $\link(F)$ has a facet $\delta$ corresponding to the face opposite $F$ in $\Delta$. If a facet contains $F$ multiple times due to self-identifications, repeat this for each occurrence of $F$.
    \item If $\Delta_1$ and $\Delta_2$ are glued in $\tri$ by a gluing which involves all vertices of $F$, then the corresponding $\delta_1$ and $\delta_2$ are glued using the same vertex identification, excluding the vertices of $F$.
  \end{itemize}
\end{definition}

This definition of a link captures the intuitive notion of a $(d-i-1)$-dimensional sphere ``around'' or ``enclosing'' $F$. If $\tri$ is a simplicial complex, it coincides with the standard definition of a link for simplicial complexes.  

The \textit{dual graph} $\Gamma(\tri)$ of $\tri$ is the graph with one node for each facet, and an arc between two nodes for each gluing between the corresponding facets. It follows that the dual graph of a generalised $d$-dimensional triangulation is a $(d+1)$-regular multigraph with loops. Moving from a triangulation to its dual graph loses all information on how the vertices are paired in each gluing, but retains some essential combinatorial structure in a form that is easy to visualise.

\subsection{Elementary Moves}
\label{subsec:moves}

In this section we describe several local modifications, or ``moves'', which can be performed on 4-manifold triangulations, modifying the triangulation without changing its underlying PL structure. These moves are the key to identifying the PL type of our triangulations, by performing a sequence of them to transform a given triangulation into a known one. The most fundamental of these moves are the so-called \textit{Pachner} or \textit{bistellar moves}:

\begin{definition}[\cite{Pachner1991}]
  Let $F$ be an internal $i$-face of a \textit{simplicial} $d$-manifold triangulation $\tri$ such that the link of $F$ is a triangulated $(n-i-1)$-sphere combinatorially isomorphic to the boundary $\partial\Delta_{n-i}$ of a standard $(n-i)$-simplex $\Delta_{n-i}$. Then the \textit{Pachner move} on $F$ replaces $\tri$ with
  \begin{align*}
    (\tri \setminus F \star \partial\Delta_{n-i}) \cup \partial F \star \Delta_{n-i}
  \end{align*}
  That is, it replaces the ``closed star'' $\overline{\openstar(F)} = F \star \partial\Delta_{n-i}$, with $\partial F \star \Delta_{n-i}$, which both have boundary $\partial F \star \partial\Delta_{n-i}$. Here $A \star B$ represents the \textit{join} of $A$ and $B$, containing, for every $p$-face of $A$ and $q$-face of $B$, a $(p+q+1)$-face spanning the vertices of those two faces.
\end{definition}

For generalised triangulations, we do not necessarily have a ``closed star'' of the form $F \star \partial\Delta_{n-i}$, since additional identifications may occur at its boundary. However, since the boundary $\partial F \star \partial\Delta_{n-i}$ does not change during a Pachner move, the move only affects $\openstar(F)$, and so we can extend this definition to generalised triangulations.
Pachner moves are fundamental even among other ``elementary'' local modifications, due to the following theorem.

\begin{theorem}[Pachner \cite{Pachner1991}]
  \label{thm:pachner}
  Two triangulations are connected by a sequence of Pachner moves if and only if they triangulate the same PL manifold.
\end{theorem}

The \textit{Pachner graph} of a PL manifold $\manifold$ is the graph with a node for each triangulation of $\manifold$ (up to combinatorial isomorphism), and an arc between two nodes if there is a Pachner move transforming the corresponding triangulations into each other. In this setting, \Cref{thm:pachner} states that the Pachner graph of $\manifold$ is connected.

In dimension 4, there are 5 types of Pachner moves, corresponding to the case where $F$ is a vertex, edge, triangle, tetrahedron or pentachoron. These are called the 5-1, 4-2, 3-3, 2-4 and 1-5 moves, respectively, corresponding to the number of pentachora in the affected region before and after the move.

\begin{figure}
\centering
\resizebox{\linewidth}{!}{\tikzset{
    >=latex,
    centremark/.style={decoration={markings, mark=at position 0.5 with {\arrow{#1}}}, postaction={decorate}}
}

\begin{tikzpicture}[every node/.style={circle, inner sep=0, outer sep=0, minimum width=0.1cm, fill=black}, every path/.style={thick}, label distance=2mm]

\node[red] (c1) at (0,0) {} ;
\path (c1) ++(90:2) node (t1) {} ;
\path (c1) ++(162:2) node (tl1) {} ;
\path (c1) ++(234:2) node (bl1) {} ;
\path (c1) ++(306:2) node (br1) {} ;
\path (c1) ++(18:2) node (tr1) {} ;

\draw (c1) -- (t1) ;
\draw (c1) -- (tl1) ;
\draw (c1) -- (bl1) ;
\draw (c1) -- (br1) ;
\draw (c1) -- (tr1) ;

\node[fill=none,outer sep=0.5] (al) at (2.5,0) {} ;
\node[fill=none,outer sep=0.5] (ar) at (3.5,0) {} ;

\draw[->] (al.90) -- (ar.90) {} ;
\draw[->] (ar.270) -- (al.270) {} ;

\node[fill=none] at (3,0.4) {1-5} ;
\node[fill=none] at (3,-0.4) {5-1} ;

\node[fill=none] (c2) at (6,0) {} ;
\path (c2) ++(90:1)  node[red] (t2) {} ;
\path (c2) ++(162:1) node[red] (tl2) {} ;
\path (c2) ++(234:1) node[red] (bl2) {} ;
\path (c2) ++(306:1) node[red] (br2) {} ;
\path (c2) ++(18:1)  node[red] (tr2) {} ;

\path (t2) ++(90:1) node (ot2) {} ;
\path (tl2) ++(162:1) node (otl2) {} ;
\path (bl2) ++(234:1) node (obl2) {} ;
\path (br2) ++(306:1) node (obr2) {} ;
\path (tr2) ++(18:1) node (otr2) {} ;

\draw (t2)  -- (ot2) ;
\draw (tl2) -- (otl2) ;
\draw (bl2) -- (obl2) ;
\draw (br2) -- (obr2) ;
\draw (tr2) -- (otr2) ;

\draw (t2) -- (tl2) ;
\draw (t2) -- (bl2) ;
\draw (t2) -- (br2) ;
\draw (t2) -- (tr2) ;

\draw (tl2) -- (bl2) ;
\draw (tl2) -- (br2) ;
\draw (tl2) -- (tr2) ;

\draw (bl2) -- (br2) ;
\draw (bl2) -- (tr2) ;

\draw (br2) -- (tr2) ;

\end{tikzpicture}

\begin{tikzpicture}[every node/.style={circle, inner sep=0, outer sep=0, minimum width=0.1cm, fill=black}, every path/.style={thick}, label distance=2mm]

\node[fill=none] (lspace) at (-2,0) {} ;
\node[fill=none] (rspace) at (8,0) {} ;

\node[red] (t1) at (0,1) {} ;
\node[red] (b1) at (0,-1) {} ;

\path (t1) ++(150:1) node (tfl1) {} ;
\path (t1) ++(110:1) node (tml1) {} ;
\path (t1) ++(70:1)  node (tmr1) {} ;
\path (t1) ++(30:1)  node (tfr1) {} ;

\path (b1) ++(-150:1) node (bfl1) {} ;
\path (b1) ++(-110:1) node (bml1) {} ;
\path (b1) ++(-70:1)  node (bmr1) {} ;
\path (b1) ++(-30:1)  node (bfr1) {} ;

\draw (t1) -- (b1) ;

\draw (t1) -- (tfl1) ;
\draw (t1) -- (tml1) ;
\draw (t1) -- (tmr1) ;
\draw (t1) -- (tfr1) ;

\draw (b1) -- (bfl1) ;
\draw (b1) -- (bml1) ;
\draw (b1) -- (bmr1) ;
\draw (b1) -- (bfr1) ;

\node[fill=none,outer sep=0.5] (al) at (2.5,0) {} ;
\node[fill=none,outer sep=0.5] (ar) at (3.5,0) {} ;

\draw[->] (al.90) -- (ar.90) {} ;
\draw[->] (ar.270) -- (al.270) {} ;

\node[fill=none] at (3,0.4) {2-4} ;
\node[fill=none] at (3,-0.4) {4-2} ;

\node[red] (bl2) at (4.5,-1) {} ;
\node[red] (tl2) at (5.5,1) {} ;
\node[red] (br2) at (6.5,-1) {} ;
\node[red] (tr2) at (7.5,1) {} ;

\path (bl2) ++(90:1) node (blu2) {} ;
\path (tl2) ++(90:1) node (tlu2) {} ;
\path (br2) ++(90:1) node (bru2) {} ;
\path (tr2) ++(90:1) node (tru2) {} ;

\path (bl2) ++(-90:1) node (bld2) {} ;
\path (tl2) ++(-90:1) node (tld2) {} ;
\path (br2) ++(-90:1) node (brd2) {} ;
\path (tr2) ++(-90:1) node (trd2) {} ;

\draw (bl2) -- (tl2) ;
\draw (bl2) -- (br2) ;
\draw (bl2) -- (tr2) ;

\draw (tl2) -- (br2) ;
\draw (tl2) -- (tr2) ;

\draw (br2) -- (tr2) ;

\draw (bl2) -- (blu2) ;
\draw (bl2) -- (bld2) ;

\draw (tl2) -- (tlu2) ;
\draw (tl2) -- (tld2) ;

\draw (br2) -- (bru2) ;
\draw (br2) -- (brd2) ;

\draw (tr2) -- (tru2) ;
\draw (tr2) -- (trd2) ;

\end{tikzpicture}}
\caption{Effect of 1-5 / 5-1 and 2-4 / 4-2 moves on the dual graph of a 4-dimensional triangulation. Red nodes show pentachora replaced by the move.}
\end{figure}
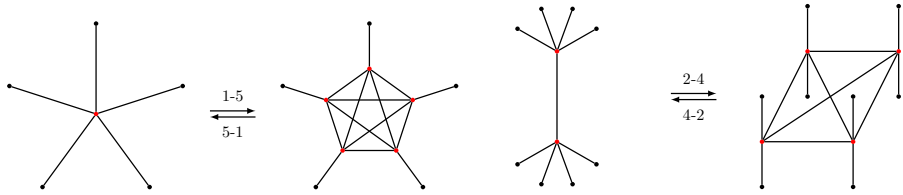

In addition to Pachner moves, there are a number of other moves we can use to modify a triangulation without changing its PL structure. There are two types of such moves on $4$-dimensional triangulations which we will make use of in the following sections:

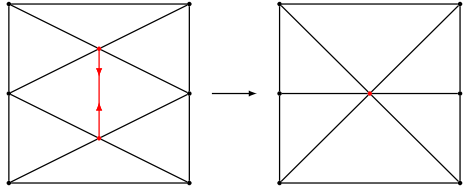
\begin{figure}
\centering
\resizebox{0.5\linewidth}{!}{\tikzset{
    >=latex,
    centremark/.style={decoration={markings, mark=at position 0.5 with {\arrow{#1}}}, postaction={decorate}},
    firstmark/.style={decoration={markings, mark=at position 0.3 with {\arrow{#1}}}, postaction={decorate}},
    secondmark/.style={decoration={markings, mark=at position 0.7 with {\arrow{#1}}}, postaction={decorate}},
}

\begin{tikzpicture}[every node/.style={circle, inner sep=0, outer sep=0, minimum width=0.1cm, fill=black}, every path/.style={thick}, label distance=2mm]

\node (BL1) at (0,0) {} ;
\node (ML1) at (0,2) {} ;
\node (TL1) at (0,4) {} ;

\node[red] (BC1) at (2,1) {} ;
\node[red] (TC1) at (2,3) {} ;

\node (BR1) at (4,0) {} ;
\node (MR1) at (4,2) {} ;
\node (TR1) at (4,4) {} ;

\draw (BL1) -- (ML1) ;
\draw (ML1) -- (TL1) ;

\draw[red,firstmark={>},secondmark={<}] (TC1) -- (BC1) ;

\draw (BR1) -- (MR1) ;
\draw (MR1) -- (TR1) ;

\draw (BL1) -- (BR1) ;
\draw (TL1) -- (TR1) ;

\draw (BL1) -- (BC1) ;
\draw (ML1) -- (BC1) ;
\draw (ML1) -- (TC1) ;
\draw (TL1) -- (TC1) ;

\draw (BR1) -- (BC1) ;
\draw (MR1) -- (BC1) ;
\draw (MR1) -- (TC1) ;
\draw (TR1) -- (TC1) ;

\draw[->] (4.5,2) -- (5.5,2) ;

\node (BL2) at (6,0) {} ;
\node (ML2) at (6,2) {} ;
\node (TL2) at (6,4) {} ;

\node[red] (C2) at (8,2) {} ;

\node (BR2) at (10,0) {} ;
\node (MR2) at (10,2) {} ;
\node (TR2) at (10,4) {} ;

\draw (BL2) -- (ML2) ;
\draw (ML2) -- (TL2) ;

\draw (BR2) -- (MR2) ;
\draw (MR2) -- (TR2) ;

\draw (BL2) -- (BR2) ;
\draw (TL2) -- (TR2) ;

\draw (BL2) -- (C2) ;
\draw (ML2) -- (C2) ;
\draw (TL2) -- (C2) ;

\draw (BR2) -- (C2) ;
\draw (MR2) -- (C2) ;
\draw (TR2) -- (C2) ;

\end{tikzpicture}}
\caption{Edge collapse for a 2-dimensional triangulation.}
\label{fig:edge-collpase}
\end{figure}

\paragraph*{Edge collapse.} Take an edge $e$ between two distinct vertices and contract it to a point, while also collapsing all of the higher dimensional faces containing $e$ along with it, so that triangles are flattened to edges, tetrahedra are flattened to triangles and so on - see \Cref{fig:edge-collpase}. This reduces the number of vertices by one while also reducing the number of $d$-simplices. Note that certain combinatorial conditions exist to ensure an edge collapse does not change the PL structure of the underlying triangulation \cite{Burton2012-Regina}.

\begin{figure}
\centering
\resizebox{0.9\linewidth}{!}{\tikzset{
    >=latex,
    centremark/.style={decoration={markings, mark=at position 0.5 with {\arrow{#1}}}, postaction={decorate}},
    mydash/.style={dash=on 1pt off 0.5pt phase 0pt}
}

\begin{tikzpicture}[every node/.style={circle, inner sep=0, outer sep=0, minimum width=0.1cm, fill=black}, every path/.style={thick}, label distance=2mm, on grid]



\node[fill=none] (c1) at (0,0) {} ;
\node (l1) [fill=none, left=of c1] {} ;
\node (r1) [fill=none, right=of c1] {} ;

\node (t1) [above=1.5 of c1] {} ;
\node (b1) [below=1.5 of c1] {} ;

\path (l1) arc (180:360:1 and 0.4) node[pos=0,auto,anchor=center] (sl1) {} node[pos=1,auto,anchor=center] (sr1) {} ;
\path (l1) arc (180:0:1 and 0.5) node[pos=0.6,auto,anchor=center] (sm1) {} ;

\draw[gray] (l1) arc (180:0:1 and 0.7) (r1) {} ;
\draw[gray] (l1) arc (180:0:1 and -0.5) (r1) {} ;

\draw (l1) -- (r1) ;
\draw (l1) -- (sm1) ;
\draw (sm1) -- (r1) ;

\path (c1) -- node[pos=0.33,auto,anchor=center,red] (v) {} (sm1);

\draw (v) -- (l1) ;
\draw (v) -- (r1) ;
\draw (v) -- (sm1);

\draw (l1) to[bend left=20] (t1) ;
\draw[dashed] (t1) to[bend left=10] (sm1) ;
\draw (r1) to[bend right=20] (t1) ;

\draw (l1) to[bend right=20] (b1) ;
\draw[dashed] (sm1) to[bend left=5] (b1) ;
\draw (r1) to[bend left=20] (b1) ;


\node[fill=none,outer sep=0.5] (al) at (1.5,0) {} ;
\node[fill=none,outer sep=0.5] (ar) at (2.5,0) {} ;

\draw[->,black] (al) -- (ar) {} ;

\node[fill=none] at (2,0.3) {\small{2-0 (vertex)}} ;


\node[fill=none] (c2) at (4,0) {} ;
\node (l2) [fill=none, left=of c2] {} ;
\node (r2) [fill=none, right=of c2] {} ;

\node (t2) [above=1.5 of c2] {} ;
\node (b2) [below=1.5 of c2] {} ;

\path (l2) arc (180:360:1 and 0.4) node[pos=0,auto,anchor=center] (sl2) {} node[pos=1,auto,anchor=center] (sr2) {} ;
\path (l2) arc (180:0:1 and 0.5) node[pos=0.6,auto,anchor=center] (sm2) {} ;

\draw (l2) -- (r2) ;
\draw (l2) -- (sm2) ;
\draw (sm2) -- (r2) ;

\draw (l2) to[bend left=20] (t2) ;
\draw[dashed] (t2) to[bend left=10] (sm2) ;
\draw (r2) to[bend right=20] (t2) ;

\draw (l2) to[bend right=20] (b2) ;
\draw[dashed] (sm2) to[bend left=5] (b2) ;
\draw (r2) to[bend left=20] (b2) ;



\node[fill=none] (c3) at (7,0) {} ;
\path (c3) ++(200:1.25) node (l3) {} ;
\path (c3) ++(20:1.25) node (r3) {} ;
\path (c3) ++(110:0.5) node (b3) {} ;
\path (c3) ++(290:0.5) node (f3) {} ;

\draw (l3) -- (f3);
\draw (f3) -- (r3);
\draw (r3) -- (b3);
\draw (b3) -- (l3);

\draw[red,mydash] (l3) -- (r3);

\draw (f3) to[out=90,in=90] (b3);
\draw[dash=on 1pt off 1pt phase 0pt] (f3) to[out=-90,in=-90] (b3);

\draw[gray] (l3) to[bend left=50] (r3) ;
\draw[gray] (l3) to[bend right=55] (r3) ;

\path (c3) ++(155:1.7) node (oul3) {} ;
\path (c3) ++(55:1.7) node (our3) {} ;
\path (c3) ++(240:1.7) node (oll3) {} ;
\path (c3) ++(-25:1.7) node (olr3) {} ;

\draw (l3) -- (oul3) ;
\draw (f3) to[bend right=20] (oul3) ;
\draw (b3) to[bend right=30] (oul3) ;

\draw (r3) -- (our3) ;
\draw (f3) to[bend left=20] (our3) ;
\draw (b3) to[bend left=40] (our3) ;

\draw (l3) -- (oll3) ;
\draw (f3) to[bend left=20] (oll3) ;
\draw[mydash] (b3) to[bend left=10] (oll3) ;

\draw (r3) -- (olr3) ;
\draw (f3) to[bend right=20] (olr3) ;
\draw[mydash] (b3) to[bend right=10] (olr3) ;


\node[fill=none,outer sep=0.5] (al) at (8.75,0) {} ;
\node[fill=none,outer sep=0.5] (ar) at (9.75,0) {} ;

\draw[->,black] (al) -- (ar) {} ;

\node[fill=none] at (9.125,0.3) {\small{2-0 (edge)}} ;


\node[fill=none] (c4) at (11.5,0) {} ;
\path (c4) ++(200:1.25) node (l4) {} ;
\path (c4) ++(20:1.25) node (r4) {} ;
\path (c4) ++(110:0.5) node (b4) {} ;
\path (c4) ++(290:0.5) node (f4) {} ;

\draw (l4) -- (f4);
\draw (f4) -- (r4);
\draw (r4) -- (b4);
\draw (b4) -- (l4);

\draw (f4) -- (b4);

\path (c4) ++(155:1.7) node (oul4) {} ;
\path (c4) ++(55:1.7) node (our4) {} ;
\path (c4) ++(240:1.7) node (oll4) {} ;
\path (c4) ++(-25:1.7) node (olr4) {} ;

\draw (l4) -- (oul4) ;
\draw (f4) to[bend right=20] (oul4) ;
\draw (b4) to[bend right=30] (oul4) ;

\draw (r4) -- (our4) ;
\draw (f4) to[bend left=20] (our4) ;
\draw (b4) to[bend left=40] (our4) ;

\draw (l4) -- (oll4) ;
\draw (f4) to[bend left=20] (oll4) ;
\draw[mydash] (b4) to[bend left=10] (oll4) ;

\draw (r4) -- (olr4) ;
\draw (f4) to[bend right=20] (olr4) ;
\draw[mydash] (b4) to[bend right=10] (olr4) ;

\end{tikzpicture}}
\caption{2-0 vertex and edge moves for a 3-dimensional triangulation.}
\label{fig:2-0-3dim}
\end{figure}
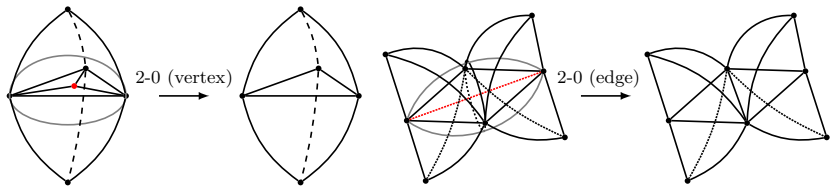

\paragraph*{2-0 vertex, edge, and triangle moves.} Each of these three moves takes two pentachora enclosing an internal face and ``flattens'' them, resulting in a triangulation with two fewer pentachora. The 2-0 vertex move takes two pentachora joined along four faces and enclosing a vertex, and pushes the two remaining ridges together to flatten them to a single tetrahedron. Similarly, the 2-0 edge (triangle) move takes two pentachora joined along three (two) faces to enclose an edge (triangle), pushing the four (six) remaining ridges down onto two (three). The corresponding moves in dimension 3 are illustrated in \Cref{fig:2-0-3dim}. Again, under certain combinatorial conditions, each of the three types of 2-0 moves indeed does not change the PL structure of the underlying space. These conditions can easily be checked, and are detailed in \cite{Burton2012-Regina} along with the documentation and source code for Regina \cite{regina}.

\section{Setup}
\label{sec:setup}

\subsection{The Connected Sums of $\CP^2$ and $S^2 \times S^2$}
\label{subsec:cp2-s2xs2-csums}

We now review implications of the results from \Cref{subsec:manifold} to classify all \textit{smoothable} simply connected closed 4-manifolds that arise as connected sums of $\CP^2$ and $S^2 \times S^2$, cf. \Cref{lem:cp2-s2xs2-toptypes}. 
The results that follow are all well-established in the literature, but often stated without proof \cite{BasakSpreer2016,SpreerTillmann2018}, so we take the time here to go through them in full detail.

First, we aim to rearrange all connected sums of $\CP^2$ and $S^2 \times S^2$ into a standard form. Aside from the standard properties of the connected sum with respect to commutativity, associativity and orientation reversal, we make use of a result due to Kirby \cite{Kirby1989}, stating that $\CP^2 \# (S^2 \times S^2)$ is PL-homeomorphic to $\CP^2 \# \CP^2 \# \overline{\CP^2}$.

\begin{lemma}
  Let $\manifold$ be a PL 4-manifold which can be constructed as a connected sum of finitely many copies of the standard PL $\CP^2$ and $S^2 \times S^2$. Then $\manifold$ is PL-homeomorphic to exactly one of the following:
  \begin{enumerate}[label=(\Roman*)]
    \item $S^4$ \label{item:standard-form-sphere}
    \item $\#^k \CP^2$ for each $k \ge 1$ \label{item:standard-form-cp2}
    \item $\#^k S^2 \times S^2$ for each $k \ge 1$ \label{item:standard-form-s2xs2}
    \item $\#^k (\CP^2 \# \overline{\CP^2})$, for each $k \ge 1$ \label{item:standard-form-cp2pm}
    \item $(\#^k \CP^2) \# (\#^l S^2 \times S^2)$, for each $k,l \ge 1$ \label{item:standard-form-cp2s2xs2}
  \end{enumerate}
  \label{lem:family-standard-forms}
\end{lemma}
\begin{proof}
  Generically, the connected sum may involve any number of copies of $\CP^2$, $\overline{\CP^2}$, $S^2 \times S^2$ and $\overline{S^2 \times S^2}$. However, for any $\manifold'$ we have $\manifold' \# (S^2 \times S^2) \cong \manifold' \# (\overline{S^2 \times S^2})$ since $S^2 \times S^2$ has an orientation-reversing self-diffeomorphism (extended from the orientation-reversing self-diffeomorphism on $S^2$) and so we may disregard $\overline{S^2 \times S^2}$. By commutativity, we have:
  \begin{align*}
    \manifold &\cong (\#^a \CP^2) \# (\#^b \overline{\CP^2}) \# (\#^c S^2 \times S^2)
  \end{align*}
  for some $a,b,c \in \NN$, where we may assume freely that $a \ge b$ by reversing the orientation of $\manifold$ if necessary. We then have three cases.

  \smallskip
  
  \noindent
  \textbf{Case 1: $a=b=0$}
  \begin{align*}
    \manifold &\cong \#^c S^2 \times S^2
  \end{align*}
  which is either \labelcref{item:standard-form-sphere} or \labelcref{item:standard-form-s2xs2}.

  \smallskip
  
  \noindent
  \textbf{Case 2: $a=b>0$}
  \begin{align*}
     \manifold &\cong (\#^{a-1} \CP^2) \# (\#^a \overline{\CP^2}) \# \CP^2 \# (\#^c S^2 \times S^2) \\
               &\cong (\#^{a-1} \CP^2) \# (\#^a \overline{\CP^2}) \# \CP^2 \# (\#^c (\CP^2 \# \overline{\CP^2})) \\
               &\cong \#^{a+c} (\CP^2 \# \overline{\CP^2})
  \end{align*}
  which is \labelcref{item:standard-form-cp2pm}.

  \smallskip
  
  \noindent
  \textbf{Case 3: $a>b$}
  \begin{align*}
    \manifold &\cong (\#^{a-b} \CP^2) \# (\#^b (\CP^2 \# \overline{\CP^2})) \# (\#^c S^2 \times S^2) \\
              &\cong (\#^{a-b} \CP^2) \# (\#^b S^2 \times S^2) \# (\#^c S^2 \times S^2) \\
              &\cong (\#^{a-b} \CP^2) \# (\#^{b+c} S^2 \times S^2)
  \end{align*}
  which is either \labelcref{item:standard-form-cp2} or \labelcref{item:standard-form-cp2s2xs2}.

Finally, each of the manifolds from \labelcref{item:standard-form-sphere}-\labelcref{item:standard-form-cp2s2xs2} are pairwise non-homeomorphic, and hence also non-PL-homeomorphic. This can be seen directly from their intersection forms using \hyperref[thm:freedman]{Freedman's theorem} -- \labelcref{item:standard-form-sphere} has an even form of rank 0 and signature 0, \labelcref{item:standard-form-cp2} has an odd form of rank $k$ and signature $k$, \labelcref{item:standard-form-s2xs2} has an even form of rank $2k$ and signature 0, \labelcref{item:standard-form-cp2pm} has an odd form of rank $2k$ and signature 0, and finally \labelcref{item:standard-form-cp2s2xs2} has an odd form of rank $k+2l$ and signature $k$ (where $0 < k < k+2l$).
\end{proof}

With this classification in hand, we turn our attention to the scope of intersection forms (and thus topological types) realisable by these connected sums:

\begin{lemma}
  Every simply connected, closed, topological 4-manifold $\manifold$ which
  \begin{enumerate}[label=(\alph*)]
    \item admits a PL-structure, and
    \item does not have an even intersection form with nonzero signature
  \end{enumerate}
  can be obtained as a connected sum of finitely many copies of the standard PL $\CP^2$ and $S^2 \times S^2$.
  \label{lem:cp2-s2xs2-toptypes}
\end{lemma}
\begin{proof}
  By \Cref{cor:rank-sig-type}, the topological type of $\manifold$ is determined the rank, signature and type of its intersection form $Q$. Hence, we only need to show that there is a manifold listed in \Cref{lem:family-standard-forms} whose intersection form has this same rank, signature and type.
  \smallskip
  
  \noindent
  \textbf{Case 1: $Q$ is even.}
  By assumption, $\sign(Q)=0$. The series \labelcref{item:standard-form-sphere} and \labelcref{item:standard-form-s2xs2} produce even intersection forms of any possible rank with signature zero.

  \noindent
  \textbf{Case 2: $Q$ is odd and indefinite.}
  By \Cref{lem:odd-indefinite-forms}, 
  \begin{align*}
    Q \cong \frac{\rank(Q)+\sign(Q)}{2} \cdot (+1) \oplus \frac{\rank(Q)-\sign(Q)}{2} \cdot (-1)
  \end{align*}
  which is the intersection form of the element of series \labelcref{item:standard-form-cp2s2xs2} with $k=\sign(Q)$ and $l = (\rank(Q)-\sign(Q))/2$.

  \noindent
  \textbf{Case 3: $Q$ is odd and definite.}
  Since $\manifold$ is only orientable (since it is simply connected) and not oriented, we may assume without loss of generality that $Q$ is positive definite. Since $Q$ is the intersection form of a smooth manifold, by \hyperref[thm:donaldson]{Donaldson's theorem} $Q \cong \rank(Q) \cdot (+1)$. This is the intersection form of \labelcref{item:standard-form-cp2} with $k=\rank(Q)$.
\end{proof}

Finally, if we extend our connected sums to also include the (PL standard) $K3$ surface, we obtain a family of $4$-manifolds often termed the \textit{standard} simply connected PL 4-manifolds. The final lemma for this section implies, modulo \hyperref[conj:11/8]{the 11/8 conjecture}, that including this one additional type of component is sufficient to extend \Cref{lem:cp2-s2xs2-toptypes} and remove the restriction on the intersection form entirely. This suggests that, if a triangulation of $K3$ with $2\beta_2(K3)+2 = 46$ were to be found, it may be possible to extend the families presented in this paper to include all topological types of smoothable simply connected closed 4-manifolds, using similar techniques to those used in \Cref{sec:series,sec:PLtype}. Attempts to find such a 46-pentachoron $K3$ triangulation have been made by Burke, who has produced the currently smallest-known $K3$ triangulation with $54$ pentachora \cite{Burke2024}.

\begin{lemma}
  If the $11/8$ Conjecture is true, then every simply connected, closed, topological 4-manifold $\manifold$ which admits a PL-structure can be obtained as a connected sum of finitely many copies of the standard PL $\CP^2$, $S^2 \times S^2$ and $K3$ surface.
  \label{lem:11-over-8-connected-sums}
\end{lemma}
\begin{proof}
  Let $Q$ be the intersection form of $\manifold$. The only remaining case not covered in \Cref{lem:cp2-s2xs2-toptypes} is when $Q$ is even and $\sign(Q) \neq 0$ (without loss of generality, $\sign(Q) > 0$). By \hyperref[thm:rokhlin]{Rokhlin's theorem}, $\sign(Q) = 16k$ for some $k \in \NN$, and hence, by \hyperref[conj:11/8]{the 11/8 conjecture}, $\rank(Q) \ge 22k$. Since the rank and signature must differ by a multiple of 2, we have $\rank(Q)=22k+2l$ for some $l \in \NN$. The even forms with this rank and signature are realised by the connected sums $(\#^k K3) \# (\#^l S^2 \times S^2)$.
\end{proof}

\subsection{Connected Sums of 4-Manifold Triangulations}
\label{subsec:csum-method}

For our proof of the PL type of $\CPplusfamily{k}$, $\CPplusminusfamily{k}$, $\SSfamily{k}$ and $\CPSSfamily{k}{l}$, we need a method to construct connected sums of 4-manifolds in the triangulation setting. This is what enables a proof by induction, allowing us, for example, to construct a certified triangulation of $\#^k\CP^2$ from already-verified triangulations of $\#^{k-1}\CP^2$ and $\CP^2$. Moreover, for the computational parts of the proof to be practical we require that our method add as few additional pentachora as possible.

For simplicial complexes this is straightforward: simply remove a pentachoron (which is an embedded 4-ball) from each triangulation and identify the resulting boundaries. However, in our generalised triangulations, the pentachora may have self-identifications, so we require a different approach. What follows is an adaptation of the method used in Regina \cite{regina} for 3-manifolds.

Abstractly, the procedure being performed is as follows. In each of our summands $\tri_1$ and $\tri_2$, we ``pop open'' a ridge, un-gluing the interiors of the two tetrahedra which had been glued together but leaving all their lower-dimensional faces identified. This has the effect of puncturing the manifold, and the new $S^3$ boundary component consists of two tetrahedra identified along all of their faces (the ``pillow'' triangulation of $S^3$). We then take a triangulation $\mathcal{C}$ of $S^3 \times I$ where both boundary components are pillow $S^3$'s, and glue one of its boundary components to the puncture in $\tri_1$ and the other to the puncture in $\tri_2$.

However, we cannot actually perform the ``pop open" operation in a triangulation as we have defined it. A triangulation is defined entirely by the gluings between tetrahedra, so when a gluing is removed all lower-dimensional identifications of their faces will be lost as well. So we must ensure that these identifications are re-introduced when we glue in $\mathcal{C}$.

We choose $\mathcal{C}$ such that every gluing uses the identity. The procedure is:
\begin{enumerate}[label=(\roman*)]
  \item Choose one gluing from $\tri_1$ -- where tetrahedron $t_1$ is glued to tetrahedron $t_2$ according to a permutation $\sigma \in \Sym_n$ -- and remove it.
  \item Let the two tetrahedra of the first boundary component of $\mathcal{C}$ be $c_1$ and $c_2$. Glue $c_1$ to $t_1$ by a permutation $\rho$, and glue $c_2$ to $t_2$ by $\sigma \rho$.
  \item Repeat steps (i)-(ii) using $\tri_2$ and the second boundary component of $\mathcal{C}$. We obtain the connected sum with opposite orientation by swapping which tetrahedron is labelled $c_1$ and which is labelled $c_2$. 
\end{enumerate}
Observe that the vertex of $t_1$ with label $i$ is identified with the vertex of $t_2$ with label $(\sigma \rho \, \id \, \rho^{-1})(i) = \sigma(i)$, the same as it was originally identified with in $\tri_1$. Hence all vertices, edges and triangles of $\tri_1$ and $\tri_2$ will be identified in the same way before and after the procedure, as required. The only requirement on the permutation $\rho$ is that it indeed maps the index of $c_1$ to the index of $t_1$ (or else they could not be glued along this permutation) -- since $\sigma$ maps the index of $t_1$ to that of $t_2$, $\sigma \rho$ will correctly map the index of $c_2$ (the same as $c_1$) to that of $t_2$. For simplicity, we choose $\rho$ to be the permutation mapping the vertices of $c_1$ in ascending order to those of $t_1$ in ascending order.

It remains to specify the triangulation $\mathcal{C}$. Triangulate two tetrahedral prisms with four pentachora each, analogously to the standard (3-dimensional) triangulation of a triangular prism with 3 tetrahedra. Then glue together the ``sides'' (but not the ``ends'') of the two prisms by the identity map. This results in the gluing table shown in \Cref{tab:cylinder}. We can verify that $\mathcal{C}$ is $S^3 \times I$ by finding a sequence of Pachner moves from it to a verified triangulation of $S^3 \times I$. The verified triangulation is produced by using Regina's \texttt{iBundle()} function to triangulate a direct sum of the pillow $S^3$ with $I$, and the full move sequence is included in \Cref{app:sequences}. 

\begin{table}
\centering
\begin{tabular}{|c|c|c|c|c|c|}
\hline
Pentachoron & Tet 0123 & Tet 0124 & Tet 0134 & Tet 0234 & Tet 1234 \\ \hline \hline
0           & -        & 1 (0124) & 4 (0134) & 4 (0234) & 4 (1234) \\ \hline
1           & 5 (0123) & 0 (0124) & 2 (0134) & 5 (0234) & 5 (1234) \\ \hline
2           & 6 (0123) & 6 (0124) & 1 (0134) & 3 (0234) & 6 (1234) \\ \hline
3           & 7 (0123) & 7 (0124) & 7 (0134) & 2 (0234) & -        \\ \hline
4           & -        & 5 (0124) & 0 (0134) & 0 (0234) & 0 (1234) \\ \hline
5           & 1 (0123) & 4 (0124) & 6 (0134) & 1 (0234) & 1 (1234) \\ \hline
6           & 2 (0123) & 2 (0124) & 5 (0134) & 7 (0234) & 2 (1234) \\ \hline
7           & 3 (0123) & 3 (0124) & 3 (0134) & 6 (0234) & -        \\ \hline
\end{tabular}
\caption{Gluing table for the triangulation $\mathcal{C}$ of $S^3 \times I$. The entry ``$p$ ($abcd$)'' indicates that this pentachoron is glued to pentachoron $p$ along a tetrahedron such that its (ordered) vertices are paired with vertices $a,b,c,d$ of $p$. Blank entries indicate a boundary tetrahedron.}
\label{tab:cylinder}
\end{table}

\section{The triangulations $\CPplusfamily{k}$, $\SSfamily{k}$, $\CPplusminusfamily{k}$ and $\CPSSfamily{k}{l}$}
\label{sec:series}

In this section we describe four infinite series of triangulated $4$-manifolds denoted by $\CPplusfamily{k}$, $\SSfamily{k}$, $\CPplusminusfamily{k}$ and $\CPSSfamily{k}{l}$ (for $k,l \in \NN$). In \Cref{sec:PLtype}, their members are shown to be PL-homeomorphic to the connected sums of $\CP^2$ and $S^2 \times S^2$ listed in cases \labelcref{item:standard-form-cp2}-\labelcref{item:standard-form-cp2s2xs2} of \Cref{lem:family-standard-forms}. That is, $\CPplusfamily{k}$, $\SSfamily{k}$, $\CPplusminusfamily{k}$ and $\CPSSfamily{k}{l}$ contain triangulations of all connected sums of $\CP^2$ and $S^2 \times S^2$ with the piecewise linear standard structure. A complete description of each of these series as gluing tables is given in \Cref{tab:s4,tab:cp2k,tab:s2xs2k,tab:cp2PMk,tab:cp2k-s2xs2l}.

The triangulations are constructed as chains of pentachora, see \Cref{fig:series-examples} for example illustrations of their dual graphs. There are three basic units, each of which has two boundary tetrahedra. The first unit has two pentachora and corresponds to a $\CP^2$ component in the connected sum, which we call a \textit{bow} $\mathcal{B}$. The second and third both have four pentachora and correspond to an $S^2 \times S^2$, and a $\CP^2 \# \overline{\CP^2}$ summand. We call these latter triangulations an \textit{even hook} $\mathcal{H}_0$ and an \textit{odd hook} $\mathcal{H}_1$ respectively. Note that both hooks have identical dual graphs, but are not combinatorially isomorphic. 

To build the overall chain, start with the $1$-pentachoron ``double snapped'' $4$-ball triangulation described in \Cref{tab:double-loop}, which we call $\textrm{DSB}_2$ (following the terminology of Burke in his thesis\footnote{note the difference in terminology from \cite{SpreerTobin24FaceNumbers}, where this was referred to as $DS_1$.}). Note that the single boundary component of $\textrm{DSB}_2$ is a $1$-tetrahedron triangulation of $S^3$ which we call $\mathcal{S}$ (see \Cref{tab:boundary}). To realise a particular connected sum $\manifold$ of $\CP^2$ and $S^2 \times S^2$ summands, take an appropriate number of copies of bows and/or hooks (referred to as {\em units} in what follows), and put them in sequence with all bows appearing first. Gluing the boundary of $\textrm{DSB}_2$ to the first boundary tetrahedron of the first unit leaves, again, a $1$-tetrahedron boundary. Continue in this way, gluing one boundary tetrahedron of the next unit onto the existing $1$-tetrahedron boundary until all units are joined up. Finally, the chain is completed by gluing a second copy of $\textrm{DSB}_2$ to the remaining $1$-tetrahedron boundary. The exact gluings involved in this process are listed in \Cref{tab:s4,tab:s2xs2k,tab:cp2k,tab:cp2k-s2xs2l}.

By construction, the resulting triangulation has $2$ pentachora for the two copies of $\textrm{DSB}_2$, plus $2=2\beta_2(\CP^2)$ for each $\CP^2$ component in the connected sum, $4=2\beta_2(S^2 \times S^2)$ for each $S^2 \times S^2$ component, and $4=2\beta_2(\CP^2 \# \overline{\CP^2})$ for each $\CP^2 \# \overline{\CP^2}$ component. Hence it has a total of $2\beta_2(\manifold)+2$ pentachora.

We conclude this section by observing that each unit, and in fact each family, has a non-trivial symmetry. Moreover, all of these symmetries restricted to the same symmetry on the $1$-tetrahedron gluing sites. This provides a choice of two gluings in each step of building up a chain. This fact will be useful in determining the piecewise linear homeomorphism types of our triangulations in \Cref{sec:PLtype}.

\begin{lemma}
  Let $\tri$ be any member triangulation of $\CPplusfamily{k}$, $\SSfamily{k}$, $\CPplusminusfamily{k}$ and $\CPSSfamily{k}{l}$, and let $\sigma$ be a gluing between two adjacent units of $\tri$. Then
  \begin{enumerate}[label=(\alph*)]
    \item removing $\sigma$ produces two connected components which both have boundaries that are combinatorially isomorphic to $\mathcal{S}$;
    \item there exists exactly one other permutation $\rho$ which, when used in place of $\sigma$, produces a combinatorially isomorphic triangulation; and
    \item using any other permutation in place of $\sigma$ or $\rho$ results in an invalid triangulation with an edge identified with itself in reverse.
  \end{enumerate}
  \label{lem:series-symmetry}
\end{lemma}
\begin{proof}
  Constructing the chain (in either direction), we claim that at every step of the construction the $1$-tetrahedron boundary is always isomorphic to $\mathcal{S}$. This is relatively straightforward to verify inductively, with the following cases to check: gluing a unit to $\textrm{DSB}_2$ at the ``front'' of the chain, gluing a unit to $\textrm{DSB}_2$ at the ``back'' of the chain, gluing a unit to itself, gluing an $S^2 \times S^2$ unit to a $\CP^2$ unit, or gluing a $\CP^2$ unit to an $S^2 \times S^2$ unit. When un-gluing two adjacent units, one component is a partial chain constructed front-to-back and the other one is a partial chain constructed back-to-front, both with boundary $\mathcal{S}$.

  In $\tri$, suppose the gluing from the preceding chain to the $i$-th unit $\mathcal{U}$ is $\sigma \in \Sym_5$. By the above, any such gluing (an element of a subgroup isomorphic to $\Sym_4$, where the tetrahedra used in the gluing are fixed) corresponds to affinely identifying two copies of $\mathcal{S}$. There are $4!=24$ ways to do this, $22$ of which result in an edge being identified with itself in reverse, and hence an invalid triangulation. One of the remaining two gluings is $\sigma$, the other we denote by $\rho$. Let $\tri'$ be the triangulation obtained from $\tri$ by replacing $\sigma$ with $\rho$.

  It is straightforward to check that each of the three units and $\textrm{DSB}_2$ have exactly one non-trivial symmetry -- that is, a relabelling of its pentachora and their vertices resulting in an \textit{identical}, not just combinatorially isomorphic, triangulation -- and that this non-identity symmetry maps the first boundary tetrahedron of the unit to itself while relabelling its vertices by a non-identity permutation. 

  We construct a combinatorial isomorphism from $\tri$ to $\tri'$ as follows. First, relabel $\tri$ by applying the above non-identity symmetry $S_\mathcal{U}$ on $\mathcal{U}$, and the identity relabelling on the rest of the triangulation. Since $S_\mathcal{U}$ is a symmetry of $\mathcal{U}$, this leaves the entire triangulation unchanged with the possible exception of the two gluings joining $\mathcal{U}$ to the units (or $\textrm{DSB}_2$) immediately before and after it. Moreover, the first of these \textit{must} change, since the preceeding chain has the same labelling while the tetrahedron it is glued to in $\mathcal{U}$ is relabelled non-trivially -- this has the effect of composing $\sigma$ with a non-identity permutation. The only other permutation possible for this gluing is $\rho$, so this relabelling must replace $\sigma$ with $\rho$.

  If the gluing joining $\mathcal{U}$ to the rest of the chain after it is unchanged, then we have obtained $\tri'$ and we are done. If not, then by the same argument as above there are exactly two possibilities for this gluing, and we can apply the same process on the $(i+1)$-st unit to swap it to the desired permutation. Again, this leaves the triangulation otherwise unchanged, with the exception of the gluing between the $(i+1)$-st and $(i+2)$-nd units. Continuing in this fashion until we reach the end of the chain (potentially using the symmetry of $\textrm{DSB}_2$ in the final step), we obtain $\tri'$ from $\tri$ by only relabelling, hence they are indeed combinatorially isomorphic. 
\end{proof}

\begin{figure}
\centering
\resizebox{0.7\linewidth}{!}{\begin{tikzpicture}[every node/.style={circle, inner sep=0, outer sep=0, minimum width=0.1cm, fill=black}, every path/.style={thick}, label distance=2mm]

  \node(0) [label=below:0] {} ;
  \node(1) [label=below:1, right=of 0] {} ;
  \node(2) [label=below:2, right=of 1] {} ;
  \node(3) [label=below:3, right=of 2] {} ;
  \node(4) [label=below:4, right=of 3] {} ;
  \node(5) [label=below:5, right=of 4] {} ;
  \node(6) [label=below:6, right=of 5] {} ;
  \node(7) [label=below:7, right=of 6] {} ;
  \node(8) [label=below:8, right=of 7] {} ;
  \node(9) [label=below:9, right=of 8] {} ;

  \draw[scale=1.5] (0) to[out=150,in=-150,loop] (0) ;
  \draw[scale=3] (0) to[out=130,in=-130,loop] (0) ;

  \draw (0) -- (1) ;

  \draw (1) to[bend left=30] (2) ;
  \draw (1) to[bend right=30] (2) ;
  \draw[scale=1.5] (1) to[out=60,in=120,loop] (1) ;
  \draw[scale=1.5] (2) to[out=60,in=120,loop] (2) ;

  \draw (2) -- (3) ;

  \draw (3) to[bend left=30] (4) ;
  \draw (3) to[bend right=30] (4) ;
  \draw[scale=1.5] (3) to[out=60,in=120,loop] (3) ;
  \draw[scale=1.5] (4) to[out=60,in=120,loop] (4) ;

  \draw (4) -- (5) ;

  \draw (5) to[bend left=30] (6) ;
  \draw (5) to[bend right=30] (6) ;
  \draw[scale=1.5] (5) to[out=60,in=120,loop] (5) ;
  \draw[scale=1.5] (6) to[out=60,in=120,loop] (6) ;

  \draw (6) -- (7) ;

  \draw (7) to[bend left=30] (8) ;
  \draw (7) to[bend right=30] (8) ;
  \draw[scale=1.5] (7) to[out=60,in=120,loop] (7) ;
  \draw[scale=1.5] (8) to[out=60,in=120,loop] (8) ;

  \draw (8) -- (9) ;

  \draw[scale=1.5] (9) to[out=30,in=-30,loop] (9) ; 
  \draw[scale=3] (9) to[out=50,in=-50,loop] (9) ;

\end{tikzpicture}}
\resizebox{0.7\linewidth}{!}{\begin{tikzpicture}[every node/.style={circle, inner sep=0, outer sep=0, minimum width=0.1cm, fill=black}, every path/.style={thick}, label distance=2mm]

  \node(0)  [label=above:0]                    {} ;
  \node(1)  [label=above:1,right=of 0]         {} ;
  \node(3)  [label=above right:3,right=of 1]   {} ;
  \node(5)  [label=above:5,right=of 3]         {} ;
  \node(7)  [label=above right:7,right=of 5]   {} ;
  \node(9)  [label=above:9,right=of 7]         {} ;
  \node(11) [label=above right:11,right=of 9]  {} ;
  \node(13) [label=above:13,right=of 11]       {} ;
  \node(15) [label=above right:15,right=of 13] {} ;
  \node(17) [label=above:17,right=of 15]       {} ;

  \node(2)  [label=below:2,below=of 1]         {} ;
  \node(4)  [label=below:4,below=of 3]         {} ;
  \node(6)  [label=below:6,below=of 5]         {} ;
  \node(8)  [label=below:8,below=of 7]         {} ;
  \node(10) [label=below:10,below=of 9]         {} ;
  \node(12) [label=below:12,below=of 11]        {} ;
  \node(14) [label=below:14,below=of 13]       {} ;
  \node(16) [label=below:16,below=of 15]       {} ;

  \draw[scale=1.5] (0) to[out=150,in=-150,loop] (0) ;
  \draw[scale=3]   (0) to[out=130,in=-130,loop] (0) ;

  \draw[scale=1.5] (2) to[out=150,in=-150,loop] (2) ;
  \draw[scale=1.5] (3) to[out=60,in=120,loop]   (3) ;
  \draw[scale=1.5] (4) to[out=60,in=120,loop]   (4) ;
  \draw[scale=3]   (4) to[out=40,in=140,loop]   (4) ;
  \draw (0) -- (1) ;
  \draw (1) to[bend left=30] (2) ;
  \draw (1) to[bend right=30] (2) ;
  \draw (2) -- (4) ;
  \draw (1) to[bend left=30] (3) ;
  \draw (1) to[bend right=30] (3) ;

  \draw[scale=1.5] (6) to[out=150,in=-150,loop] (6) ;
  \draw[scale=1.5] (7) to[out=60,in=120,loop]   (7) ;
  \draw[scale=1.5] (8) to[out=60,in=120,loop]   (8) ;
  \draw[scale=3]   (8) to[out=40,in=140,loop]   (8) ;
  \draw (3) -- (5) ;
  \draw (5) to[bend left=30] (6) ;
  \draw (5) to[bend right=30] (6) ;
  \draw (6) -- (8) ;
  \draw (5) to[bend left=30] (7) ;
  \draw (5) to[bend right=30] (7) ;

  \draw[scale=1.5] (10) to[out=150,in=-150,loop] (10) ;
  \draw[scale=1.5] (11) to[out=60,in=120,loop]   (11) ;
  \draw[scale=1.5] (12) to[out=60,in=120,loop]   (12) ;
  \draw[scale=3]   (12) to[out=40,in=140,loop]   (12) ;
  \draw (7) -- (9) ;
  \draw (9) to[bend left=30] (10) ;
  \draw (9) to[bend right=30] (10) ;
  \draw (10) -- (12) ;
  \draw (9) to[bend left=30] (11) ;
  \draw (9) to[bend right=30] (11) ;

  \draw[scale=1.5] (14) to[out=150,in=-150,loop] (14) ;
  \draw[scale=1.5] (15) to[out=60,in=120,loop]   (15) ;
  \draw[scale=1.5] (16) to[out=60,in=120,loop]   (16) ;
  \draw[scale=3]   (16) to[out=40,in=140,loop]   (16) ;
  \draw (11) -- (13) ;
  \draw (13) to[bend left=30] (14) ;
  \draw (13) to[bend right=30] (14) ;
  \draw (14) -- (16) ;
  \draw (13) to[bend left=30] (15) ;
  \draw (13) to[bend right=30] (15) ;

  \draw[scale=1.5] (17) to[out=30,in=-30,loop] (17) ;
  \draw[scale=3]   (17) to[out=50,in=-50,loop] (17) ;
  \draw (15) -- (17) ;

\end{tikzpicture}}
\resizebox{0.7\linewidth}{!}{\begin{tikzpicture}[every node/.style={circle, inner sep=0, outer sep=0, minimum width=0.1cm, fill=black}, every path/.style={thick}, label distance=2mm]

  \node(0)  [label=above:0]                   {} ;
  \node(1)  [label=below:1,right=of 0]        {} ;
  \node(2)  [label=below:2,right=of 1]        {} ;
  \node(3)  [label=below:3,right=of 2]        {} ;
  \node(4)  [label=below:4,right=of 3]        {} ;
  \node(5)  [label=above:5,right=of 4]        {} ;
  \node(7)  [label=above right:7,right=of 5]  {} ;
  \node(9)  [label=above:9,right=of 7]        {} ;
  \node(11) [label=above right:11,right=of 9] {} ;
  \node(13) [label=above:13,right=of 11]      {} ;
  
  \node(6)  [label=below:6,below=of 5]        {} ;
  \node(8)  [label=below:8,below=of 7]        {} ;
  \node(10) [label=below:10,below=of 9]       {} ;
  \node(12) [label=below:12,below=of 11]      {} ;

  \draw[scale=1.5] (0) to[out=150,in=-150,loop] (0) ;
  \draw[scale=3]   (0) to[out=130,in=-130,loop] (0) ;

  \draw[scale=1.5] (1) to[out=60,in=120,loop] (1) ;
  \draw[scale=1.5] (2) to[out=60,in=120,loop] (2) ;
  \draw (0) -- (1) ;
  \draw (1) to[bend left=30]  (2) ;
  \draw (1) to[bend right=30] (2) ;

  \draw[scale=1.5] (3) to[out=60,in=120,loop] (3) ;
  \draw[scale=1.5] (4) to[out=60,in=120,loop] (4) ;
  \draw (2) -- (3) ;
  \draw (3) to[bend left=30]  (4) ;
  \draw (3) to[bend right=30] (4) ;

  \draw[scale=1.5] (6) to[out=150,in=-150,loop] (6) ;
  \draw[scale=1.5] (7) to[out=60,in=120,loop]   (7) ;
  \draw[scale=1.5] (8) to[out=60,in=120,loop]   (8) ;
  \draw[scale=3]   (8) to[out=40,in=140,loop]   (8) ;
  \draw (4) -- (5) ;
  \draw (5) to[bend left=30]  (6) ;
  \draw (5) to[bend right=30] (6) ;
  \draw (6) -- (8) ;
  \draw (5) to[bend left=30]  (7) ;
  \draw (5) to[bend right=30] (7) ;

  \draw[scale=1.5] (10) to[out=150,in=-150,loop] (10) ;
  \draw[scale=1.5] (11) to[out=60,in=120,loop]   (11) ;
  \draw[scale=1.5] (12) to[out=60,in=120,loop]   (12) ;
  \draw[scale=3]   (12) to[out=40,in=140,loop]   (12) ;
  \draw (7) -- (9) ;
  \draw (9) to[bend left=30]  (10) ;
  \draw (9) to[bend right=30] (10) ;
  \draw (10) -- (12) ;
  \draw (9) to[bend left=30]  (11) ;
  \draw (9) to[bend right=30] (11) ;
  
  \draw[scale=1.5] (13) to[out=30,in=-30,loop] (13) ;
  \draw[scale=3]   (13) to[out=50,in=-50,loop] (13) ;
  \draw (11) -- (13) ;

\end{tikzpicture}}
\caption{Dual graphs of the triangulations $\CPplusfamily{4}$ (top), $\CPplusminusfamily{4}$ and $\SSfamily{4}$ (middle), and $\CPSSfamily{2}{2}$ (bottom).}
\label{fig:series-examples}
\end{figure}

\begin{table}
\centering
\begin{tabular}{|c|r|r|r|r|r|}
\hline
Pentachoron & Tet 0123 & Tet 0124 & Tet 0134 & Tet 0234 & Tet 1234 \\ \hline \hline
0           & -        & 0 (0324) & 0 (3214) & 0 (0214) & 0 (3104) \\ \hline
\end{tabular}
\caption{Gluing table for $\textrm{DSB}_2$.}
\label{tab:double-loop}
\end{table}

\begin{table}
\centering
\begin{tabular}{|c|r|r|r|r|r|}
\hline
Tetrahedron & Tri 012 & Tri 013 & Tri 023 & Tri 123 \\ \hline \hline
0           & 0 (032) & 0 (321) & 0 (021) & 0 (310) \\ \hline
\end{tabular}
\caption{Gluing table for $\mathcal{S}$.}
\label{tab:boundary}
\end{table}

\begin{table}
\centering
\begin{tabular}{|c|c|c|c|c|c|}
\hline
Pentachoron & Tet 0123 & Tet 0124 & Tet 0134 & Tet 0234 & Tet 1234 \\ \hline \hline
0           & 1 (0123) & 0 (0324) & 0 (3214) & 0 (0214) & 0 (3104) \\ \hline
0           & 0 (0123) & 1 (0324) & 1 (3214) & 1 (0214) & 1 (3104) \\ \hline
\end{tabular}
\caption{Gluing table for $\CPplusfamily{0}=\SSfamily{0}=\CPplusminusfamily{0}=\CPSSfamily{0}{0}$.}
\label{tab:s4}
\end{table}

\begin{table}
\centering
\begin{tabular}{|c|c|c|c|c|c|}
\hline
Pentachoron & Tet 0123 & Tet 0124 & Tet 0134 & Tet 0234 & Tet 1234 \\ \hline \hline
0           & 1 (0123) & 0 (0324) & 0 (3214) & 0 (0214) & 0 (3104) \\ \hline \hline
1           & 0 (0123) & 2 (1034) & 1 (3214) & 2 (0234) & 1 (3104) \\ \hline
2           & 3 (1203) & 2 (3124) & 1 (1024) & 1 (0234) & 2 (1204) \\ \hline \hline
3           & 2 (2013) & 4 (1034) & 3 (3214) & 4 (0234) & 3 (3104) \\ \hline
4           & 5 (1203) & 4 (3124) & 3 (1024) & 3 (0234) & 4 (1204) \\ \hline \hline
\vdots      & \vdots   & \vdots   & \vdots   & \vdots   & \vdots   \\ \hline \hline
\multirow{2}{*}{$2k-1$} & $2k-2$ & $2k$   & $2k-1$ & $2k$   & $2k-1$ \\ 
                        & (2013) & (1034) & (3214) & (0234) & (3104) \\ \hline
\multirow{2}{*}{$2k$}   & $2k+1$ & $2k$   & $2k-1$ & $2k-1$ & $2k$   \\ 
                        & (0123) & (3124) & (1024) & (0234) & (1204) \\ \hline \hline
\multirow{2}{*}{$2k+1$} & $2k$   & $2k+1$ & $2k+1$ & $2k+1$ & $2k+1$ \\ 
                        & (0123) & (3124) & (3024) & (1304) & (1204) \\ \hline
\end{tabular}
\caption{Gluing table for $\CPplusfamily{k}$.}
\label{tab:cp2k}
\end{table}

\begin{table}
\centering
\begin{tabular}{|c|c|c|c|c|c|}
\hline
Pentachoron & Tet 0123  & Tet 0124  & Tet 0134  & Tet 0234  & Tet 1234  \\ \hline \hline
 0          &  1 (0123) &  0 (0324) &  0 (3214) &  0 (0214) &  0 (3104) \\ \hline \hline
 1          &  0 (0123) &  3 (0124) &  2 (3204) &  3 (0234) &  2 (1234) \\ \hline
 2          &  4 (0123) &  2 (1304) &  2 (2014) &  1 (3104) &  1 (1234) \\ \hline
 3          &  5 (0123) &  1 (0124) &  3 (3214) &  1 (0234) &  3 (3104) \\ \hline
 4          &  2 (0123) &  4 (1304) &  4 (2014) &  4 (1234) &  4 (0234) \\ \hline \hline
 5          &  3 (0123) &  7 (0124) &  6 (3204) &  7 (0234) &  6 (1234) \\ \hline
 6          &  8 (0123) &  6 (1304) &  6 (2014) &  5 (3104) &  5 (1234) \\ \hline
 7          &  9 (0123) &  5 (0124) &  7 (3214) &  5 (0234) &  7 (3104) \\ \hline
 8          &  6 (0123) &  8 (1304) &  8 (2014) &  8 (1234) &  8 (0234) \\ \hline \hline
\vdots      & \vdots    & \vdots    & \vdots    & \vdots    & \vdots    \\ \hline \hline
\multirow{2}{*}{$4k-3$}  & $4k-5$ & $4k-1$ & $4k-2$ & $4k-1$ & $4k-2$ \\ 
                         & (0123) & (0124) & (3204) & (0234) & (1234) \\ \hline
\multirow{2}{*}{$4k-2$}  & $4k$   & $4k-2$ & $4k-2$ & $4k-3$ & $4k-3$ \\
                         & (0123) & (1304) & (2014) & (3104) & (1234) \\ \hline
\multirow{2}{*}{$4k-1$}  & $4k+1$ & $4k-3$ & $4k-1$ & $4k-3$ & $4k-1$ \\
                         & (0123) & (0124) & (3214) & (0234) & (3104) \\ \hline
\multirow{2}{*}{$4k$  }  & $4k-2$ & $4k$   & $4k$   & $4k$   & $4k$   \\
                         & (0123) & (1304) & (2014) & (1234) & (0234) \\ \hline \hline
\multirow{2}{*}{$4k+1$}  & $4k-1$ & $4k+1$ & $4k+1$ & $4k+1$ & $4k+1$ \\
                         & (0123) & (0324) & (3214) & (0214) & (3104) \\ \hline
\end{tabular}
\caption{Gluing table for $\CPplusminusfamily{k}$.}
\label{tab:cp2PMk}
\end{table}

\begin{table}
\centering
\begin{tabular}{|c|c|c|c|c|c|}
\hline
Pentachoron & Tet 0123  & Tet 0124  & Tet 0134  & Tet 0234  & Tet 1234  \\ \hline \hline
 0          &  1 (0123) &  0 (0324) &  0 (3214) &  0 (0214) &  0 (3104) \\ \hline \hline

 1          &  0 (0123) &  3 (0124) &  2 (3204) &  3 (0234) &  2 (1234) \\ \hline
 2          &  4 (0123) &  2 (3014) &  2 (1204) &  1 (3104) &  1 (1234) \\ \hline
 3          &  5 (0123) &  1 (0124) &  3 (3214) &  1 (0234) &  3 (3104) \\ \hline
 4          &  2 (0123) &  4 (3014) &  4 (1204) &  4 (1234) &  4 (0234) \\ \hline \hline

 5          &  3 (0123) &  7 (0124) &  6 (3204) &  7 (0234) &  6 (1234) \\ \hline
 6          &  8 (0123) &  6 (3014) &  6 (1204) &  5 (3104) &  5 (1234) \\ \hline
 7          &  9 (0123) &  5 (0124) &  7 (3214) &  5 (0234) &  7 (3104) \\ \hline
 8          &  6 (0123) &  8 (3014) &  8 (1204) &  8 (1234) &  8 (0234) \\ \hline \hline

\vdots      & \vdots    & \vdots    & \vdots    & \vdots    & \vdots    \\ \hline \hline

\multirow{2}{*}{$4k-3$}  & $4k-5$ & $4k-1$ & $4k-2$ & $4k-1$ & $4k-2$ \\ 
                         & (0123) & (0124) & (3204) & (0234) & (1234) \\ \hline
\multirow{2}{*}{$4k-2$}  & $4k$   & $4k-2$ & $4k-2$ & $4k-3$ & $4k-3$ \\ 
                         & (0123) & (3014) & (1204) & (3104) & (1234) \\ \hline
\multirow{2}{*}{$4k-1$}  & $4k+1$ & $4k-3$ & $4k-1$ & $4k-3$ & $4k-1$ \\ 
                         & (0123) & (0124) & (3214) & (0234) & (3104) \\ \hline
\multirow{2}{*}{$4k$  }  & $4k-2$ & $4k$   & $4k$   & $4k$   & $4k$   \\ 
                         & (0123) & (3014) & (1204) & (1234) & (0234) \\ \hline \hline

\multirow{2}{*}{$4k+1$}  & $4k-1$ & $4k+1$ & $4k+1$ & $4k+1$ & $4k+1$ \\ 
                         & (0123) & (0324) & (3214) & (0214) & (3104) \\ \hline
\end{tabular}
\caption{Gluing table for $\SSfamily{k}$.}
\label{tab:s2xs2k}
\end{table}

\begin{table}
\centering
\begin{tabular}{|c|c|c|c|c|c|}
\hline
Pentachoron & Tet 0123  & Tet 0124  & Tet 0134  & Tet 0234  & Tet 1234  \\ \hline \hline

0          &  1 (0123) &  0 (0324) &  0 (3214) &  0 (0214) &  0 (3104) \\ \hline \hline

1          &  0 (0123) &  2 (1034) &  1 (3214) &  2 (0234) &  1 (3104) \\ \hline
2          &  3 (1203) &  2 (3124) &  1 (1024) &  1 (0234) &  2 (1204) \\ \hline \hline

 \vdots     & \vdots    & \vdots    & \vdots    & \vdots    & \vdots    \\ \hline \hline

\multirow{2}{*}{$2k-3$} & $2k-4$ & $2k-2$ & $2k-3$ & $2k-2$ & $2k-3$ \\
                        & (2013) & (1034) & (3214) & (0234) & (3104) \\ \hline
\multirow{2}{*}{$2k-2$} & $2k-1$ & $2k-2$ & $2k-3$ & $2k-3$ & $2k-2$ \\ 
                        & (1203) & (3124) & (1024) & (0234) & (1204) \\ \hline \hline

\multirow{2}{*}{$2k-1$} & $2k-2$ & $2k$   & $2k-1$ & $2k$   & $2k-1$ \\
                        & (2013) & (1034) & (3214) & (0234) & (3104) \\ \hline
\multirow{2}{*}{$2k$  } & $2k+1$ & $2k$   & $2k-1$ & $2k-1$ & $2k$   \\
                        & (1204) & (3124) & (1024) & (0234) & (1204) \\ \hline
\multirow{2}{*}{$2k+1$} & $2k+3$ & $2k$   & $2k+2$ & $2k+3$ & $2k+2$ \\
                        & (0123) & (2013) & (4230) & (0234) & (1234) \\ \hline
\multirow{2}{*}{$2k+2$} & $2k+2$ & $2k+4$ & $2k+2$ & $2k+1$ & $2k+1$ \\
                        & (4013) & (0124) & (1230) & (4130) & (1234) \\ \hline
\multirow{2}{*}{$2k+3$} & $2k+1$ & $2k+5$ & $2k+3$ & $2k+1$ & $2k+3$ \\
                        & (0123) & (0124) & (4231) & (0234) & (4130) \\ \hline
\multirow{2}{*}{$2k+4$} & $2k+4$ & $2k+2$ & $2k+4$ & $2k+4$ & $2k+4$ \\ 
                        & (4013) & (0124) & (1230) & (1234) & (0234) \\ \hline \hline

\multirow{2}{*}{$2k+5$} & $2k+7$ & $2k+3$ & $2k+6$ & $2k+7$ & $2k+6$ \\
                        & (0123) & (0124) & (4230) & (0234) & (1234) \\ \hline
\multirow{2}{*}{$2k+6$} & $2k+6$ & $2k+8$ & $2k+6$ & $2k+5$ & $2k+5$ \\
                        & (4013) & (0124) & (1230) & (4130) & (1234) \\ \hline
\multirow{2}{*}{$2k+7$} & $2k+5$ & $2k+9$ & $2k+7$ & $2k+5$ & $2k+7$ \\
                        & (0123) & (0124) & (4231) & (0234) & (4130) \\ \hline
\multirow{2}{*}{$2k+8$} & $2k+8$ & $2k+6$ & $2k+8$ & $2k+8$ & $2k+8$ \\ 
                        & (4013) & (0124) & (1230) & (1234) & (0234) \\ \hline \hline

\vdots     & \vdots    & \vdots    & \vdots    & \vdots    & \vdots    \\ \hline \hline

\multirow{2}{*}{$2k+4l-3$} & $2k+4l-1$ & $2k+4l-5$ & $2k+4l-2$ & $2k+4l-1$ & $2k+4l-2$ \\ 
                           & (0123)    & (0124)    & (4230)    & (0234)    & (1234)    \\ \hline
\multirow{2}{*}{$2k+4l-2$} & $2k+4l-2$ & $2k+4l$   & $2k+4l-2$ & $2k+4l-3$ & $2k+4l-3$ \\ 
                           & (4013)    & (0124)    & (1230)    & (4130)    & (1234)    \\ \hline
\multirow{2}{*}{$2k+4l-1$} & $2k+4l-3$ & $2k+4l+1$ & $2k+4l-1$ & $2k+4l-3$ & $2k+4l-1$ \\ 
                           & (0123)    & (0124)    & (4231)    & (0234)    & (4130)    \\ \hline
\multirow{2}{*}{$2k+4l$  } & $2k+4l$   & $2k+4l-2$ & $2k+4l$   & $2k+4l$   & $2k+4l$   \\ 
                           & (4013)    & (0124)    & (1230)    & (1234)    & (0234)    \\ \hline \hline

\multirow{2}{*}{$2k+4l+1$} & $2k+4l+1$ & $2k+4l-1$ & $2k+4l+1$ & $2k+4l+1$ & $2k+4l+1$ \\ 
                           & (0423)    & (0124)    & (4231)    & (0231)    & (4130)    \\ \hline
\end{tabular}
\caption{Gluing table for $\CPSSfamily{k}{l}$.}
\label{tab:cp2k-s2xs2l}
\end{table}

\section{Identifying the PL Types of $\CPplusfamily{k}$, $\SSfamily{k}$, $\CPplusminusfamily{k}$ and $\CPSSfamily{k}{l}$}
\label{sec:PLtype}

In this section we prove that the triangulations constructed in \Cref{sec:series} are in one-to-one correspondence with the manifolds described in \Cref{lem:family-standard-forms,lem:cp2-s2xs2-toptypes}. We start with the following statement.

\begin{theorem}
  Let $\CP^2$ and $S^2 \times S^2$ be endowed with the standard PL structure. The following PL-homeomorphisms hold: 
  \begin{enumerate}[label=(\Roman*)]
    \item $\CPplusfamily{0} = \SSfamily{0} = \CPplusminusfamily{0} = \CPSSfamily{0}{0} \cong S^4$ \label{item:series-sphere}
    \item $\CPplusfamily{k} \cong \#^k \CP^2$ for all $k \ge 1$ \label{item:series-cp2}
    \item $\SSfamily{k} \cong \#^k S^2 \times S^2$ for all $k \ge 1$ \label{item:series-s2xs2}
    \item $\CPplusminusfamily{k} \cong \#^k (\CP^2 \# \overline{\CP^2})$, for all $k \ge 1$ \label{item:series-cp2pm}
    \item $\CPSSfamily{k}{l} \cong (\#^k \CP^2) \# (\#^l S^2 \times S^2)$, for all pairs $k,l \ge 1$ \label{item:series-cp2s2xs2}
  \end{enumerate}
  \label{thm:main}
\end{theorem}

\begin{corollary}
  \label{cor:pltype}
  Every simply connected, closed, topological 4-manifold $\manifold$ which:
  \begin{enumerate}[label=(\alph*)]
    \item admits a PL-structure, and
    \item does not have an even intersection form with nonzero signature
  \end{enumerate}
  is homeomorphic to exactly one of the triangulations $\CPplusfamily{k}$ ($k \ge 0$), $\SSfamily{k}$ ($k \ge 1$), $\CPplusminusfamily{k}$ ($k \ge 1$) or $\CPSSfamily{k}{l}$ ($k,l \ge 1$). In particular, $\manifold$ admits a triangulation with $2\beta_2(\manifold)+2$ pentachora.
\end{corollary}
\begin{proof}[Proof of \Cref{cor:pltype} and \Cref{thm:introduction}]
  Direct corollary of \Cref{thm:main} and \Cref{lem:family-standard-forms,lem:cp2-s2xs2-toptypes}.
\end{proof}

Note that, due to \Cref{cor:rank-sig-type} and once we have established simple connectivity, the topological homeomorphism type for any of our triangulations can be determined through calculating rank, signature and type of its intersection form. \Cref{thm:pachner}, on the other hand, provides us with an avenue to prove \textit{PL-homeomorphism}: if we can find a sequence of Pachner moves connecting two triangulations then we know they are PL-homeomorphic. However, searching the Pachner graph to find such a sequence becomes rapidly intractable as the number of pentachora increases -- even in a single example. Instead, our approach to establish PL-homeomorphisms for each member of our infinite families of triangulations is to reduce the problem to a finite set of PL-homeomorphisms between small triangulations, which we can then verify by computer proof.

The key is to use a single sequence of Pachner moves which works for every element of a series -- computationally it will be much more convenient to think of this as ``cutting out'' a small section of a triangulation, operating on it with Pachner moves, and gluing it back in. This relies on the simple but crucial observation that ``Pachner moves commute with gluings'', or more precisely:

\begin{lemma}
  Let $\tri$ be a triangulation of a $d$-manifold, $\tri'$ a triangulation obtained from $\tri$ by removing a gluing, and $\sigma$ a face of $\tri'$ which is a valid site for a Pachner move. Then:
  \begin{enumerate}[label=(\alph*)]
      \item \label{item:pachner-gluing-valid} $\sigma$ is also a valid site for a Pachner move in $\tri$.
      \item \label{item:pachner-gluing-bdry} Performing the Pachner move on $\sigma$ in $\tri'$ results in a triangulation $\tri''$ whose boundary is combinatorially isomorphic to that of $\tri'$.
      \item \label{item:pachner-gluing-commute} Re-adding the same gluing to $\tri''$ (in a way consistent with the original labelling of $\tri'$) yields the same result up to combinatorial isomorphism as performing the Pachner move on $\sigma$ in $\tri$.
  \end{enumerate}
\end{lemma}
\begin{proof}
  Since the Pachner move on $\sigma$ is valid, it is an internal face of $\tri'$. Moreover, all of the faces in $\openstar(\sigma)$ are also internal, since they all have $\sigma$ as a face. Adding an additional gluing only affects boundary faces, and so does not affect elements of $\openstar(\sigma)$ or how they are glued together. Hence $\link(\sigma)$ is identical in $\tri$ and $\tri'$, proving \labelcref{item:pachner-gluing-valid}.

  Now the Pachner move on $\sigma$ only affects $\openstar(\sigma)$, and thus only affects internal faces. Hence it does not affect the boundary of $\tri'$, proving \labelcref{item:pachner-gluing-bdry}. This also shows that \labelcref{item:pachner-gluing-commute} is well posed, and the result follows from the fact that the Pachner move only affects $\openstar(\sigma)$ in $\tri$ or $\tri'$, and adding the gluing only affects boundary faces in $\tri'$ or $\tri''$.
\end{proof}

\begin{corollary}
  \label{cor:unglue-reglue}
  Let $\tri_1$ be a triangulation of a $d$-manifold, and $\tri_2$ the triangulation obtained from $\tri_1$ by un-gluing any number of gluings, performing any number of Pachner moves and then re-gluing the removed gluings in the same way (with respect to the boundaries of the triangulations involved). Then there exists a sequence of Pachner moves from $\tri_1$ to $\tri_2$. In particular, they are PL-homeomorphic.
\end{corollary}

\begin{proof}[Proof of \Cref{thm:main}]

We now proceed with a proof by induction. To prove parts \labelcref{item:series-sphere}-\labelcref{item:series-cp2s2xs2} of \Cref{thm:main} simultaneously, it is sufficient to show the following base cases and induction steps:

\begin{enumerate}[label=(B\arabic*)]
  \item $\CPplusfamily{0} = \SSfamily{0} = \CPplusminusfamily{0} = \CPSSfamily{0}{0} \cong S^4$
  \item $\CPplusfamily{1} \cong \CP^2$
  \item $\SSfamily{1} \cong S^2 \times S^2$
  \item $\CPplusminusfamily{1} \cong \CP^2 \# \overline{\CP^2}$
  \item $\CPSSfamily{1}{1} \cong \CP^2 \# S^2 \times S^2$
\end{enumerate}

\begin{enumerate}[label=(I\arabic*)]
  \item adding a bow between the first $\textrm{DSB}_2$ and the first bow is PL-homeomorphic to taking a connected sum with $\CPplusfamily{1} \cong \CP^2$ \label{item:induction-cp2}
  \item adding an even hook between the last even hook and the end $\textrm{DSB}_2$ is PL-homeomorphic to taking a connected sum with $\SSfamily{1} \cong S^2 \times S^2$ \label{item:induction-s2xs2}
  \item adding an odd hook between the last odd hook and the end $\textrm{DSB}_2$ is PL-homeomorphic to taking a connected sum with $\CPplusminusfamily{1} \cong \CP^2 \# \overline{\CP^2}$ \label{item:induction-cp2pm}
\end{enumerate}

Noting that since bows always appear before (even) hooks in the order, \labelcref{item:induction-cp2} and \labelcref{item:induction-s2xs2} together prove \labelcref{item:series-cp2s2xs2}.

The base cases are a straightforward case of computer search through the Pachner graph, to join each to a known triangulation of the given PL-manifold. For $S^4$ we use the standard ``pillow'' $S^4$, for $\CP^2$, $S^2 \times S^2$ and $\CP^2 \# \overline{\CP^2}$ we use the minimal triangulations provided in Regina \cite{regina} (these are due to Burke, obtained by simplifying triangulations constructed from the corresponding Kirby diagrams \cite{Burke2024,regina}). For $\CP^2 \# S^2 \times S^2$, we take a connected sum of these triangulations by the method described in \Cref{subsec:csum-method}. For the precise triangulations and Pachner move sequences used, see \Cref{app:sequences}.

Let us first consider \labelcref{item:induction-cp2}. Let $\mathcal{T}(k)$ denote either $\CPplusfamily{k}$ or $\CPSSfamily{k}{l}$ for some given $l > 1$. Take a connected sum between $\CPplusfamily{1}$ and $\mathcal{T}(k)$ using the method of \Cref{subsec:csum-method}, choosing face 1234 of pentachoron 0 in $\CPplusfamily{1}$ and face 0134 of pentachoron 0 in $\mathcal{T}(k)$ as the site for the connected sum (``$t_1$'' in the notation of \Cref{subsec:csum-method}) and using the reversed orientation -- call this verified $\CP^2 \# \mathcal{T}(k)$ triangulation $\mathcal{V}(k+1)$. A priori it is not clear this corresponds to the ``positive'' connected sum with $\CP^2$, but since there are only two possibilities this is easily verified by simply calculating the intersection form of one such $\mathcal{V}(k+1)$. We need to show $\mathcal{T}(k+1) \cong \mathcal{V}(k+1)$.

In $\mathcal{T}(k+1)$, un-glue the gluing after the second bow to obtain two connected components -- a ``head'' $\mathcal{H}_\mathcal{T}$ consisting of one $\textrm{DSB}_2$ and two bows, and a long ``tail'' consisting of the rest of the units and the final $\textrm{DSB}_2$. In $\mathcal{V}(k+1)$, un-glue the gluing that originally came after the first bow in $\mathcal{T}(k)$, obtaining a 15-pentachoron head $\mathcal{H}_\mathcal{V}$ which contains the site of the connected sum, and a tail identical to the tail in $\mathcal{T}(k+1)$. The triangulations $\mathcal{H}_\mathcal{T}$ and $\mathcal{H}_\mathcal{V}$ do not depend on the choice of $k$ or $\mathcal{T}_k$, so we can complete the induction step by finding a sequence of elementary moves from $\mathcal{H}_\mathcal{V}$ to $\mathcal{H}_\mathcal{T}$ and applying \Cref{cor:unglue-reglue} to show $\mathcal{T}(k+1) \cong \mathcal{V}(k+1)$. An explicit sequence of 64 moves is included in \Cref{app:sequences}.

The only caveat of this procedure is that we must ensure that when we re-glue $\mathcal{H}_\mathcal{T}$ to the tail, we do this in the same way as the original gluing between $\mathcal{H}_\mathcal{V}$ and the tail. A priori this is not trivial, since when performing the computer search to construct a Pachner sequence, we cannot practically keep track of how the site for the re-gluing is relabelled. However, \Cref{lem:series-symmetry} ensures there are only two possibilities for how to perform this re-gluing, and both give combinatorially isomorphic triangulations, so in fact there is nothing to check.

The remaining cases \labelcref{item:induction-s2xs2} and \labelcref{item:induction-cp2pm} are proven in exactly the same way, with the exception that the site for the connected sum and the ``head'' occur at the end of the chain instead of the beginning. 

For these cases, the connected sum $\mathcal{V}(l+1)$ is obtained using face 1234 of the last pentachoron in $\SSfamily{l}$, $\CPplusminusfamily{l}$ or $\CPSSfamily{k}{l}$ and face 0234 of pentachoron 5 in $\SSfamily{1}$ or $\CPplusminusfamily{1}$, without reversing orientation. The gluing removed to obtain the head and tail is the one before the second-last hook in $\mathcal{T}(l+1)$ (and the matching gluing in $\mathcal{V}(l+1)$). For each of \labelcref{item:induction-s2xs2} and \labelcref{item:induction-cp2pm} we find a sequence of 80 moves, both included in \Cref{app:sequences}.

\end{proof}

\begin{figure}
\centering
\resizebox{\linewidth}{!}{\input{figures/cp2k-proof-diagram-final}}
\caption{Schematic for applying \Cref{cor:unglue-reglue} to prove \Cref{thm:main}, demonstrated for the case $\CPplusfamily{1}\#\CPSSfamily{2}{1} \cong \CPSSfamily{3}{1}$. $\mathcal{H}_\mathcal{V}$ and $\mathcal{H}_\mathcal{T}$ are the components on the left in the third and fourth diagram, respectively, and the shared tail is the component on the right in both diagrams. Arcs coloured red indicate the location of the connected sum.}
\label{fig:PL-proof-schematic}
\end{figure}

\section{Computational Methods}
\label{sec:computation}

The proof of \Cref{thm:main} in \Cref{sec:PLtype} reduces the problem of finding infinitely many PL-homeomorphisms to a finite set of computations. In this final section, we turn our attention to the practical problem of performing these computations, detailing our main algorithm under realistic time and memory restrictions.

Our computational problems are all of the same type: given two known triangulations which we claim are PL-homeomorphic, find a sequence of Pachner moves joining them. This is similar to the well-studied problem of simplifying a given triangulation \cite{Burton2012-Regina,Burton2011SoCG,Burke2024}. Na\"ively, if the triangulations are indeed PL-homeomorphic, then it is certainly possible to find a sequence by brute force, performing a breadth-first-search of the Pachner graph of one triangulation until we encounter the other. From a theoretical perspective we know that this method will eventually succeed, and indeed will always find the shortest possible such sequence. However, we have no bound on how long the sequence may be or how many triangulations will be visited before it is found. 

The possibility of an arbitrarily large shortest sequence is not a problem we can practically solve here. Instead, we maximise our chances of finding an existant connecting sequence by implementing some measures to make our search more efficient. For this, note that from a practical point of view, the bottleneck when computing connecting sequences of Pachner moves is not in time, but in space. Attempting to record every triangulation visited by a search algorithm will typically overwhelm available memory well before computation time becomes an issue.

We refer to our algorithm as ``outside-in'', developed to complete the computer proof in \Cref{thm:main}. At its core, it is a modification of the na\"ive breadth-first-search method described above. The idea is to perform two breadth-first searches simultaneously, exploring the Pachner graph from both starting triangulations and searching for an overlap between the two components (constructing the sequence from their respective ends (the ``outside'') in to the middle (``in'')). At any given step of the algorithm, we minimise the number of triangulations which must be held in memory, to be checked against newly encountered ones. By constructing paths from the starting triangulations one \emph{ring}, i.e., level of depth of the breadth-first search tree, at a time, and carefully keeping track of what has already been checked, we can check new triangulations against only the most recent two rings while still ensuring there are no duplicates in any earlier rings.

The sequences provided in \Cref{app:sequences} were produced by a combination of this algorithm, specifically the extended version with simultaneous simplification described in \Cref{subsec:outside-in-simplify}, and the ``up-side-down'' simplification algorithm developed by Burke \cite{Burke2024}. Where possible, the entire sequence was constructed with ```outside-in'' (which by design typically results in a shorter sequence), and where this proved too taxing on memory requirements, sequences were found by first using ``up-side-down'' simplifications to join the larger triangulation to one with fewer pentachora, and then using ``outside-in'' to join this smaller triangulation to the target one.

\subsection{The Basic Algorithm}
\label{subsec:outside-in-nosimplify}

We begin with a simpler version of the algorithm, which is less practical but comes with better theoretical guarantees. Note that in all of the following, triangulations are considered equal if they are combinatorially isomorphic.

\begin{algorithm_thm}[Outside-In Without Simplification]
\label{alg:outside-in-nosimplify} {\color{white} .}

\medskip

\noindent
\textbf{Input}: non-combinatorially isomorphic triangulations $\tri_1$ and $\tri_2$, $k \in \NN$.

\medskip

\noindent
\textbf{Output}: a shortest sequence of Pachner moves connecting $\tri_1$ to $\tri_2$ with all triangulations in the sequence having at most $n_{\max}+k$ facets, $n_{\max} = \max \{f_4(\tri_1), f_4(\tri_2)\}$, \texttt{false} if such a sequence does not exist.

\begin{enumerate}
  \item Set $\mathcal{L}_0 = \mathcal{R}_0 := \emptyset$, $\mathcal{L}_1 := \{\tri_1\}$, $\mathcal{R}_1 := \{\tri_2\}$, and $i = j := 1$.
  \item For each triangulation $\tri'$, $f_4(\tri') \leq n_{\max}+k$, obtained from a triangulation in $\mathcal{L}_i$ by a Pachner move:
  \begin{enumerate}
    \item If $\tri' \not \in \mathcal{L}_{i-1} \cup \mathcal{L}_{i}$, add it to the set $\mathcal{L}_{i+1}$.
    \item If $\tri' \in \mathcal{R}_{j}$, set $i:=i+1$ and go to Step 4.
  \end{enumerate}
  Go to Step 3.
  \item[2'.] For each triangulation $\tri'$, $f_4(\tri') \leq n_{\max}+k$, obtained from a triangulation in $\mathcal{R}_j$ by a Pachner move:
  \begin{enumerate}
    \item If $\tri' \not \in \mathcal{R}_{j-1} \cup \mathcal{R}_{j}$, add it to the set $\mathcal{R}_{j+1}$.
    \item If $\tri' \in \mathcal{L}_{i}$, set $j:=j+1$ and go to Step 4.
  \end{enumerate}
  Go to Step 3'.

  \item If $\mathcal{L}_{i+1}=\emptyset$: return \texttt{false}.

  If $\mathcal{L}_{i+1}\ne\emptyset$: set $i:=i+1$. If $|\mathcal{L}_i| \le |\mathcal{R}_j|$, go to Step 2. Otherwise, go to Step 2'.

  \item[3'.] If $\mathcal{R}_{j+1}=\emptyset$: return \texttt{false}.

  If $\mathcal{R}_{j+1}\ne\emptyset$: set $j:=j+1$. If $|\mathcal{L}_i| \le |\mathcal{R}_j|$, go to Step 2. Otherwise, go to Step 2'.

  \item Label the overlap triangulation from Step 2(b)/2'(b) as $\hat{\tri} \in \mathcal{L}_{i} \cap \mathcal{R}_{j}$.
  \item Starting from $\hat{\tri}$, trace back triangulations through the rings $\mathcal{L}_{a}$, $i \geq a \geq 1$, and $\mathcal{R}_{b}$, $j \geq b\geq 1$, and return the resulting sequence of Pachner moves.  
\end{enumerate}
\end{algorithm_thm}

Steps 1-3 are the core of the algorithm and confirm the existence of a sequence, while Steps 4-5 are simply a backtrace to reconstruct the explicit sequence of Pachner moves. 
This backtrace can be eliminated entirely by simply recording each triangulation alongside its ``parent'' when stored in the sets $\mathcal{L}_i$ or $\mathcal{R}_i$. However, in practice the time consumed by steps 4-5 (even if implemented na\"ively) is insignificant compared to the bottleneck of Steps 1-3, and storing parents requires valuable memory. 

Note that to efficiently check for combinatorial isomorphism of triangulations, triangulations are stored as \textit{isomorphism signatures} similar to those introduced in \cite{Burton2011SoCG}. Specifically, we use the ``ridge degrees'' isomorphism signature format implemented in Regina \cite{regina}.

\begin{lemma}
  \label{lem:outside-in-guarantees}
  At all steps in \Cref{alg:outside-in-nosimplify}:
  \begin{enumerate}[label=(\alph*)]
    \item $\mathcal{L}_a \cap \mathcal{L}_b = \emptyset$ for all $a \neq b$ for which they are defined \label{item:rings-disjoint-left}
    \item $\mathcal{R}_a \cap \mathcal{R}_b = \emptyset$ for all $a \neq b$ for which they are defined \label{item:rings-disjoint-right}
    \item $\mathcal{L}_a \cap \mathcal{R}_b =
      \begin{cases}
        \{\hat{\tri}\} \textrm{ if } a=i,b=j \textrm{ are maximal such that they are defined,} \\
        \qquad \, \textrm{and Step 4 has been reached} \\
        \emptyset \textrm{ otherwise}
      \end{cases}$ \label{item:rings-disjoint-cross}
  \end{enumerate}
\end{lemma}
\begin{proof}

  For part \labelcref{item:rings-disjoint-left}, we prove by induction on $a$ that the result holds for all $b < a$ (which proves the result without loss of generality).
  
  If $a \le 1$ the result is immediate, so assume $a \ge 2$. Let $\tri \in \mathcal{L}_a$ and suppose by contradiction there exists some $b < a$ such that $\tri \in \mathcal{L}_b$. By construction, $\tri \notin \mathcal{L}_{a-1} \cup \mathcal{L}_{a-2}$, so $b \le a-3$. Now we know that there must be a triangulation $\tri' \in \mathcal{L}_{a-1}$ such that $\tri$ and $\tri'$ are connected by a Pachner move. But since $\tri \in \mathcal{L}_b$, $\tri'$ must either be in $\mathcal{L}_{b-1}$, or it would have been found when constructing $\mathcal{L}_{b+1}$, in which case $\tri' \in \mathcal{L}_{b+1}$. This contradicts the induction hypothesis, since $\tri' \in \mathcal{L}_{a-1} \cup \mathcal{L}_{b-1}$ or $\tri' \in \mathcal{L}_{a-1} \cup \mathcal{L}_{b+1}$ where $b-1 < b+1 \leq a-2 < a-1$ and $a-1 < a$, completing the proof of \labelcref{item:rings-disjoint-left}. Part \labelcref{item:rings-disjoint-right} is symmetric and proven in exactly the same way.

  For part \labelcref{item:rings-disjoint-cross}, the first case is true by construction, so we prove the remaining case. For the base of the induction, Step 1 ensures that $\mathcal{L}_a \cap \mathcal{R}_b = \emptyset$ for $a,b \in \{0,1\}$.

  At any given step of the algorithm, let $\mathcal{L}_a$ and $\mathcal{R}_b$ be the most recently added ring on each side, where $a,b \in \NN$ with at least one strictly greater than 1. Without loss of generality, assume $\mathcal{L}_a$ was added more recently than $\mathcal{R}_b$. Suppose the result holds for all $a',b'$ where $a' < a$ and $b' \le b$. We need to prove that $\mathcal{L}_a \cap \mathcal{R}_{b'} = \emptyset$ for all $b' \le b$. Since we have not reached Step 4, Step 2(b)/2'(b) guarantees that $\mathcal{L}_a \cap \mathcal{R}_b = \emptyset$. Suppose there is some $\tri \in \mathcal{L}_a$ such that we also have $\tri \in \mathcal{R}_{a-r}$, $r > 0$. There exists some $\tri' \in \mathcal{L}_{a-1}$ which is connected to $\tri$ by a Pachner move. Since $\tri \in \mathcal{R}_{b-r}$, we have $\tri' \in \mathcal{R}_{b-r+1}$ or $\tri' \in \mathcal{R}_{b-r-1}$. Hence $\mathcal{L}_{a-1} \cap \mathcal{R}_{b-r+1} \ne \emptyset$ or $\mathcal{L}_{a-1} \cap \mathcal{R}_{b-r-1} \ne \emptyset$, contradicting the induction hypothesis since $a-1 < a$ and $b-r-1 < b-r+1 \le b$. This completes the proof of \labelcref{item:rings-disjoint-cross}.

\end{proof}

\begin{lemma}
  \Cref{alg:outside-in-nosimplify} terminates and is correct. That is, if a sequence is found it is of shortest possible length subject to the given conditions, and if \texttt{false} is returned then no such sequence exists.
\end{lemma}
\begin{proof}
  We claim that $\mathcal{L}_a$ is precisely the set of triangulations such that the shortest sequence of Pachner moves joining it to $\tri_1$, subject to no intermediate triangulation having more than $n_{\max} + k$ facets, has length $a-1$.

  Let $a \in \NN$ and $\tri \in \mathcal{L}_a$. By construction, $\tri$ can be obtained from $\tri_1$ by a sequence of $a-1$ Pachner moves. If they were joined by a shorter sequence, subject to the given conditions, then $\tri \in \mathcal{L}_{a'}$ with $a' < a$, a contradiction with \Cref{lem:outside-in-guarantees}. Hence the claim is true. By an identical argument, the symmetric claim holds for $\mathcal{R}_b$ and $\tri_2$.

  Now consider the case where a sequence is successfully found. Let $i$ and $j$ be the indices of the ``outermost'' $\mathcal{L}_a$ and $\mathcal{R}_b$, such that $\hat{\tri} \in \mathcal{L}_i \cap \mathcal{R}_j$. By \Cref{lem:outside-in-guarantees}\labelcref{item:rings-disjoint-cross}, any triangulation $\hat{\tri}'$ connected to $\tri_1$ by a sequence of at most $i-1$ moves (subject to given conditions) requires a sequence of at least $j-1$ moves (subject to given conditions) to connect it to $\tri_2$. Symmetrically, any $\hat{\tri}''$ connected to $\tri_2$ by a sequence of at most $j-1$ such moves requires such a sequence of at least $i-1$ such moves to connect it to $\tri_1$. Hence, no sequence between $\tri_1$ and $\tri_2$ (subject to given conditions) can have length less than $i+j-2$, which is the length of the sequence constructed by \Cref{alg:outside-in-nosimplify}.

  If \texttt{false} is returned, then the algorithm reached a state where either $\mathcal{L}_{i+1}=\emptyset$ or $\mathcal{R}_{j+1}=\emptyset$, without loss of generality assume $\mathcal{L}_{i+1}=\emptyset$. That is, all Pachner neighbours (subject to given conditions) of triangulations in $\mathcal{L}_i$ are in $\mathcal{L}:=\bigcup_{a=1}^i \mathcal{L}_a$, and hence over the course of the algorithm we have ensured that every neighbour of every $\mathcal{T} \in \mathcal{L}$ is also in $\mathcal{L}$. Thus, $\mathcal{L}$ is the set of all triangulations reachable from $\tri_1$ by any number of Pachner moves subject to the given conditions, and \Cref{lem:outside-in-guarantees}\labelcref{item:rings-disjoint-cross} shows that $\tri_2 \notin \mathcal{L}$, hence no such sequence from $\tri_1$ to $\tri_2$ can exist.
  
  Finally, all triangulations considered in the algorithm have at most $n_{\max}+k$ facets. Since there are only finitely many such triangulations, and by \Cref{lem:outside-in-guarantees} the algorithm considers new triangulations in every step -- or terminates, the algorithm must terminate in finite time.\footnote{In fact, we can even bound the running time of the algorithm by bounding the overall number of triangulations. But such a bound would be severely impractical.}
\end{proof}

\subsection{The Algorithm with Simplification}
\label{subsec:outside-in-simplify}

The running time and space requirements of \Cref{alg:outside-in-nosimplify}, while difficult to predict accurately, can be expected to increase rapidly with the number of facets of the two triangulations. In particular, larger triangulations require more computational resources to compute their isomorphism signatures, and tend to have more neighbours in the Pachner graph.

We can hence improve practical running times by first simplifying $\tri_1$ and $\tri_2$ (for example using ``up-down-simplify'', or exhaustive enumeration) before running \Cref{alg:outside-in-nosimplify}. Furthermore, since \Cref{alg:outside-in-nosimplify} is already performing a breadth-first search of the Pachner graph of each starting triangulation, we can simplify $\tri_1$ and $\tri_2$ while running \Cref{alg:outside-in-nosimplify} simultaneously.

Instead of restarting the algorithm whenever a simplification is found, we only reset the side on which the simplification occurred, and keep working with the already-computed rings on the other side. This measure yields another improvement in running times in practice, but comes with a loss of theoretical guarantees: the modified algorithm may ``miss'' shorter sequence between the two end triangulations, which would have been found when starting from scratch. 
We propose the following procedure:

\begin{algorithm_thm}[Outside-In With Simplification]
\label{alg:outside-in-simplify}{\color{white} .}

\noindent
\textbf{Input}: non-combinatorially isomorphic triangulations $\tri_1$ and $\tri_2$, $k \in \NN$.

\medskip

\noindent
\textbf{Output}: a sequence of Pachner moves connecting $\tri_1$ to $\tri_2$, with all triangulations in the sequence having at most $n_{\max}+k$ facets, $n_{\max} = \max \{f_4(\tri_1), f_4(\tri_2)\}$, \texttt{false} if such a sequence could not be found.

\begin{enumerate}
  \item Set $\mathcal{L}^1_0 = \mathcal{R}^1_0 := \emptyset$, $\mathcal{L}^1_1 := \{\tri_1\}$, $\mathcal{R}^1_1 := \{\tri_2\}$, and $i = j = l = m := 1$.
  \item For each triangulation $\tri'$, $f_4(\tri') \leq n_{\max}+k$, obtained from a triangulation in $\mathcal{L}^l_i$ by a Pachner move:
  \begin{enumerate}
    \item If $\tri' \not \in \mathcal{L}^l_{i-1} \cup \mathcal{L}^l_{i}$, add it to the set $\mathcal{L}^l_{i+1}$.
    \item If $\tri' \in \mathcal{R}^m_{j}$, set $i:=i+1$ and go to Step 4.
    \item If $f_4(\tri') <  f_4(\tri_1)$, set $l:=l+1$, $\tri_1 := \tri'$, $\mathcal{L}^l_0 := \emptyset$, $\mathcal{L}^l_1 := \{\tri'\}$, $n_{\max} =\max \{ f_4(\tri_1), f_4(\tri_2)\}$, and go to Step 3.
  \end{enumerate}
  Go to Step 3.
  \item[2'.] For each triangulation $\tri'$, $f_4(\tri') \leq n_{\max}+k$, obtained from a triangulation in $\mathcal{R}^m_j$ by a Pachner move:
  \begin{enumerate}
    \item If $\tri' \not \in \mathcal{R}^m_{j-1} \cup \mathcal{R}^m_{j}$, add it to the set $\mathcal{R}^m_{j+1}$.
    \item If $\tri' \in \mathcal{L}^l_{i}$, set $j:=j+1$ and go to Step 4.
    \item If $f_4(\tri') <  f_4(\tri_2)$, set $m:=m+1$, $\tri_2 := \tri'$, $\mathcal{R}^m_0 := \emptyset$, $\mathcal{R}^m_1 := \{\tri'\}$, $n_{\max} =\max \{ f_4(\tri_1), f_4(\tri_2)\}$, and go to Step 3'.
  \end{enumerate}
  Go to Step 3'.
  \item If $\mathcal{L}^l_{i+1}=\emptyset$: flag side $\mathcal{L}$ as \texttt{complete}. If $\mathcal{R}$ is also \texttt{complete}, return \texttt{false}. Otherwise, go to Step 2'.

  If $\mathcal{L}^l_{i+1}\ne\emptyset$: set $i:=i+1$. If $|\mathcal{L}^l_i| \le |\mathcal{R}^m_j|$ or $\mathcal{R}$ is \texttt{complete}, go to Step 2. Otherwise, go to Step 2'.

  \item[3'.] If $\mathcal{R}^m_{j+1}=\emptyset$: flag side $\mathcal{R}$ as \texttt{complete}. If $\mathcal{L}$ is also \texttt{complete}, return \texttt{false}. Otherwise, go to Step 2.

  If $\mathcal{R}^m_{j+1}\ne\emptyset$: set $j:=j+1$. If $|\mathcal{L}^l_i| \le |\mathcal{R}^m_j|$ and $\mathcal{L}$ is not \texttt{complete}, go to Step 2. Otherwise, go to Step 2'.

  \item Label the overlap triangulation from Step 2(b)/2'(b) as $\hat{\tri} \in \mathcal{L}^l_i \cap \mathcal{R}^m_j$.
  \item Starting from $\hat{\tri}$, trace back triangulations through the rings $\mathcal{L}^b_a$ and $\mathcal{R}^d_c$ for $l \geq b \geq 1$, $m \geq d \geq 1$, and each $a,c \ge 1$ for which they are defined, and return the resulting sequence of Pachner moves.  
\end{enumerate}

\end{algorithm_thm}

Note that \Cref{alg:outside-in-simplify} may fail in Step 3, reaching a dead end even when a sequence subject to the conditions of \Cref{alg:outside-in-nosimplify} exists. Moreover, the ``global'' guarantees of \Cref{lem:outside-in-guarantees} are replaced by guarantees on the ``local'' level between simplifications. Very roughly speaking, we can still expect the output sequence to be relatively short, even if it is no longer the shortest possible.

\begin{lemma}
  \label{lem:outside-in-guarantees-simplify}
  At all steps in \Cref{alg:outside-in-simplify}:
  \begin{enumerate}[label=(\alph*)]
    \item $\mathcal{L}^c_a \cap \mathcal{L}^c_b = \emptyset$ for all $a \ne b$ and $1 \le c \le l$.
    \item $\mathcal{R}^c_a \cap \mathcal{R}^c_b = \emptyset$ for all $a \ne b$ and $1 \le c \le m$.
  \end{enumerate}
\end{lemma}
\begin{proof}
  For part \labelcref{item:rings-disjoint-left}, we repeat exactly the argument from the proof of \Cref{lem:outside-in-guarantees}\labelcref{item:rings-disjoint-left} for each choice of $c$. For $c>1$, the base case $a,b \in \{0,1\}$ is guaranteed by Step 2(c), for $c=1$ it is guaranteed by Step 1 as before. Once again, part \labelcref{item:rings-disjoint-right} is completely symmetric.
\end{proof}

\begin{corollary}
  \Cref{alg:outside-in-simplify} terminates.
\end{corollary}
\begin{proof}
  The algorithm terminates when an overlap triangulation $\hat{\tri}$ is found, or both sides are flagged \texttt{complete}, that is when an empty ring is encountered on each side. Since the former case must happen before the latter, it suffices to analyse the latter case. \Cref{lem:outside-in-guarantees-simplify} guarantees that for any fixed $b$, no triangulation appears twice amongst the sets $\mathcal{L}^b_a$. No triangulation considered has more than $n_{\max}+k$ facets (using the original value of $n_{\max}$ at the start of the algorithm), so there are only finitely many triangulations which can appear. Hence for each $b$ there can be only finitely many nonempty $\mathcal{L}^b_a$. But there can also only be finitely many $b$, since $\tri_1$ has a finite number of facets which decreases by an integer every time $b$ increments, and cannot decrease below two. Hence there are only finitely many $\mathcal{L}^b_a$, and an empty one must be encountered in finite time. The symmetric argument for the $\mathcal{R}^d_c$ completes the proof.
\end{proof}

\bibliographystyle{plain}
\bibliography{references_lucy}

\newpage

\begin{appendices}
\crefalias{section}{appendix}

\section{Sequences of Elementary Moves}
\label{app:sequences}

We list explicit sequences of elementary moves from one triangulation to another. Moves are notated by a pair, where:
\begin{itemize}
  \item C,$i$ indicates the edge collapse on edge $i$
  \item E,$i$ indicates the 2-0 edge move on edge $i$
  \item T,$i$ indicates the 2-0 triangle move on triangle $i$
  \item $d$,$i$ with $d \in \{0,1,2,3,4\}$ indicates the Pachner move on $d$-face $i$
\end{itemize}
For example, the entry $2$,$5$ means perform a 3-3 move on triangle 5. We also include the isomorphism signatures of the triangulations in the standard \texttt{isoSig()} format used in Regina \cite{regina}.

At every step of the sequence, the labelling of faces used is the ``canonical'' labelling of a combinatorial isomorphism class used by Regina - produced by the Regina method \texttt{makeCanonical()}. To reproduce the sequence, the triangulation must be relabelled to this canonical form after each move.

\subsection*{Known $S^3 \times I$ to ``$\mathcal{C}$" from \Cref{subsec:csum-method}}

\begin{itemize}
  \item Known $S^3 \times I$ iso sig (regina: \texttt{Example4.iBundle(Example3.sphere())}):

    \texttt{-cKcfvvvvvvwLLzMvPvAPvvPwPvvAwAPMAyzPzQMzPPLLPLQwQLzwwAzAQLz \\ wwALPMMAQAMMzQQPQzAQwQLPQMQAMMAzQMzQQzAQLPALQMQPPAQQQPzQMQQQ \\ LAQMQQLAQMQLQMQAQAQgQajavanazaraAaFayaoaIasaJaPaSayaUauaYa1a \\ ya3aFaGaBa6aHaDa7aPaYaOa+aQa-aEaXaabZabbKagbMahbSa1aRakbTalb \\ Na0amb2anbRasbVaubwbWayb3aAb4aBb5a0aEb4aGbIb5aKbMbNbgb+ahbab \\ ibcbQbjbebRbSbsbubbbEbGbtbTbvbUbfbXbFbVbHbWbYbkbmbobZbqb0b1b \\ wbybnbIbKbxb2bzb3brb6bJb4bLb5b7bwbxb8bAbCb+b9b-bBbDbbcacMbNb \\ ccObecdcPbIbJbfcMbObhcgcicNbPbkcjclcncmcZbTb0bVboc1boc8b+bWb \\ qcfchcrc9bqcpcgcrcpc2b4bscsc-bbc5bucickcvcacuctcjcvctc-bwccc \\ wcxcycdcyclczcmczcncxcicAclcAcBcCcmcCcDcDcBcscEcwcEcAcEcFcGc \\ ycGcCcGcHcFcHcDcHcIcFcIcIcJcJcJcyaaaqbaaPbaaaaaaoaaavaaaaaaa \\ aaaabaaaaaaaaayayaoaaayavaaaqbPbqbaaqbaabaPbaaPbaaqbaaPbaaoa \\ vaoaaaoaaabavaaavaaayaaaPbaaaavaaabaaabaaabayaaaqbaaaaoaaaaa \\ aayaoayavayaoaaayavaaaaaqbqbbaPbPbqbaaqbaabaaaPbaaPbaaaaqbPb \\ qbaaPbaaaaoaoabavavaoaaaoaaabaaavaaavaaaaayayaaaPbPbaaaaaava \\ vaaaaababaaabaaaaabayayaaaqbqbaaaaaaoaoaaaaaaaaaaayaoayavaaa \\ yaaaqbqbbaaaPbPbaaqbaaaaPbaaaaqbPbaaaaoaoabaaavavaaaoaaaaava \\ aaaayaaaPbaaaaaavaaabaaabaaabaaayaaaqbaaaaaaoaaaaaaaaayaaaqb \\ aaPbaaaaaaoaaavaaaaaaaaabaaaaaaaaaaaaaaaaa}

  \item $\mathcal{C}$ iso sig:

    \texttt{iHgLAAMIabbdddfffhhhaaaaaaaaaaaaaaaaaaaaaa}

  \item Move sequence:

    C,0 C,119 C,33 C,1 C,8 C,77 C,1 C,64 C,1 C,58 C,49 C,38 C,9 0,8 0,9 C,24 2,26 2,18 E,9 1,19 E,9
  
\end{itemize}

\subsection*{Known $S^4$ to $\CPplusfamily{0}=\SSfamily{0}=\CPplusminusfamily{0}=\CPSSfamily{0}{0}$}

\begin{itemize}
  \item Known $S^4$ iso sig (``pillow'' $S^4$, two pentachora with boundaries identified -- regina: \texttt{Example4.sphere()}):

    \texttt{cPkbbbbaaaaaaaa}

  \item $\CPplusfamily{0}$ iso sig:

    \texttt{cAkaabb+aoa+aoa}

  \item Move sequence:
  
    3,0 2,4 3,4 0,4 3,0 2,3 3,0 1,3 2,7 2,4 3,4 0,3 3,0 2,3 1,0

\end{itemize}

\subsection*{Known $\CP^2$ to $\CPplusfamily{1}$}

\begin{itemize}
  \item Known $\CP^2$ iso sig (regina: \texttt{Example4.cp2()}):

    \texttt{eALQcaabcdddya1a1avaJaya2a}

  \item $\CPplusfamily{1}$ iso sig:

    \texttt{eAMMcaabccdd+aoa+aAaqbqbGa}

  \item Move sequence:
  
    3,7 2,7 1,4
  
\end{itemize}

\subsection*{Known $\CP^2 \# \overline{\CP^2}$ to $\CPplusminusfamily{1}$}

\begin{itemize}
  \item Known $\CP^2 \# \overline{\CP^2}$ iso sig (regina: \texttt{Example4.s2xs2Twisted()}):

    \texttt{gAMMPPaabccdeeff+aoa+aAaqbqbAa+a+aoa}

  \item $\CPplusminusfamily{1}$ iso sig:

    \texttt{gALAMQaacdcdeeff+aoawbaaYa+ayaYa+aoa}

  \item Move sequence:
  
    3,6 2,7 1,3 3,0 3,4 3,5 2,12 1,6 3,0 1,6 2,17 2,7 3,8 2,13 2,8 2,19 1,6 3,6 2,15 2,19 1,8 3,7 2,15 2,7 2,14 2,9 2,15 2,19 2,8 2,14 2,23 2,13 1,7 1,7 3,7 3,18 2,14 2,17 1,6 2,13 1,6 2,6 1,5
  
\end{itemize}

\subsection*{Known $S^2 \times S^2$ to $\SSfamily{1}$}

\begin{itemize}
  \item Known $S^2 \times S^2$ iso sig (regina: \texttt{Example4.s2xs2()}):

    \texttt{gALAMQaacdcdeeffPbgaVbaaJafbyaJafbva}

  \item $\SSfamily{1}$ iso sig:

    \texttt{gALAMQaacdcdeeff+aoawbaaGa+ayaGa+aoa}

  \item Move sequence:
  
    4,4 1,4 2,11 3,17 2,17 1,9 3,7 2,7 2,15 1,7 3,0 2,13 2,12 1,6 3,12 2,12 1,5 2,13 3,12 0,3
  
\end{itemize}

\subsection*{Known $\CP^2 \# (S^2 \times S^2)$ to $\CPSSfamily{1}{1}$}

\begin{itemize}
  \item Known $\CP^2 \# (S^2 \times S^2)$ iso sig ($\CPplusfamily{1} \# \SSfamily{1}$, using face 0 of pentachoron 0 for both):

    \texttt{sAMLMLAAPMMPMMPPaabcddeffhhhjjjlllnnnoopqqrrHbgaHb \\ aaDaRaya2aRafafafafafafafafafafafafadbgaTbDa2a2aRa}

  \item $\CPSSfamily{1}{1}$ iso sig:

    \texttt{iAMLMPAkaabcddegffghh1ava1aDaaaJaPbDaPbJafbfbva}

  \item Move sequence:
  
    C,8 C,7 C,7 3,14 2,2 3,11 3,11 2,31 2,35 T,12 3,10 2,6 2,2 2,20 2,11 2,14 2,0 2,35 2,20 2,2 2,2 2,9 2,6 2,32 2,27 2,24 2,8 2,8 2,0 2,11 2,14 2,0 2,37 E,3 2,28 E,3 E,4 2,2 3,5 2,2 T,19 2,12 T,2 3,0 3,17 2,23 1,7 3,23 2,27 2,25 2,19 2,27 2,24 2,24 2,26 2,26 2,17 2,25 2,11 1,8 3,6 2,4 2,23 2,20 2,16 1,5 3,23 1,5 2,15 3,24 2,19 1,6 3,22 2,14 2,9 2,7 2,16 1,6 3,18 2,18 2,10 2,17 2,25 2,27 2,16 2,20 2,27 2,29 2,26 2,22 2,28 2,13 1,9 2,23 2,12 1,5 1,7
  
\end{itemize}

\subsection*{$\CPplusfamily{1} \# \CPplusfamily{k}$ / $\CPplusfamily{1} \# \CPSSfamily{k}{l}$ head to $\CPplusfamily{k+1}$ / $\CPSSfamily{k+1}{l}$ head}

\begin{itemize}
  \item $\CPplusfamily{1} \# \CPplusfamily{k}$ / $\CPplusfamily{1} \# \CPSSfamily{k}{l}$ head iso sig: 

    \texttt{pAMMLPPwPMMMMiaabccdfffhhhjjjlllmmnoo2a+a2 \\ awbYayaYaYaYaYaYaYaYaYaYaYaYaYa3afanaxaba}

  \item $\CPplusfamily{k+1}$ / $\CPSSfamily{k+1}{l}$ head iso sig:

    \texttt{fAMMPaaabccdee1ava1aVbfaJadbPb}

  \item Move sequence:
  
    C,6 C,8 C,7 3,10 1,4 3,10 3,6 1,5 2,7 2,15 2,7 3,19 2,18 2,10 2,21 2,13 2,9 2,20 2,17 1,5 3,14 1,3 3,20 2,7 1,3 1,3 2,13 1,3 3,14 2,12 3,6 1,5 3,11 2,18 2,16 1,3 3,7 2,7 1,4 2,12 3,17 2,17 2,22 2,10 2,14 2,14 2,16 2,11 2,18 2,11 2,20 2,8 2,9 1,3 2,13 2,16 2,16 2,17 2,9 2,18 2,16 2,8 1,4 1,5
  
\end{itemize}

\subsection*{$\CPplusminusfamily{k} \# \CPplusminusfamily{1}$ head to $\CPplusminusfamily{k+1}$ head}

\begin{itemize}
  \item $\CPplusminusfamily{k} \# \CPplusminusfamily{1}$ head iso sig: 

    \texttt{tALAMLAAAwQPMMzIPaadccdeffhhhjjjlllnnnoopqrrssPbNaza \\ aaNaFbPbFbvaNaNaNaNaNaNaNaNaNaNaNaNakbSaSaaaEbNaPbNa}

  \item $\CPplusminusfamily{k+1}$ head iso sig:

    \texttt{jAMHwMMPkaabcddfgfghhii1ava1aVbaaJaVbaaRaJa2aRaJafa}

  \item Move sequence:

    C,9 C,8 C,9 3,11 1,5 2,10 3,16 3,19 2,26 2,16 1,8 3,7 2,12 2,15 2,19 2,18 2,21 2,27 2,8 2,15 2,19 2,18 1,6 3,0 2,27 1,7 3,27 2,20 2,26 2,27 2,16 2,21 2,21 2,25 2,22 2,27 2,23 2,29 2,28 1,7 3,18 1,9 2,14 3,19 2,20 2,28 2,10 2,15 2,14 2,18 2,26 2,16 2,27 1,8 2,14 2,19 2,12 2,24 2,15 2,16 2,10 2,24 3,23 2,27 2,26 1,8 1,7 3,30 2,18 3,19 2,21 1,4 3,15 2,29 2,12 1,5 2,10 1,9 1,5 1,5
  
\end{itemize}

\subsection*{$\SSfamily{l} \# \SSfamily{1}$ / $\CPSSfamily{k}{l} \# \SSfamily{1}$ head to $\SSfamily{l+1}$ / $\CPSSfamily{k}{l+1}$ head}

\begin{itemize}
  \item $\SSfamily{l} \# \SSfamily{1}$ / $\CPSSfamily{k}{l} \# \SSfamily{1}$ head iso sig: 

    \texttt{tALAMALAAMAPPMzIPaacdcdeefhhhjjjlllnnnoopqrrssYaoaAa \\aaGaYaqbGaYaoaoaoaoaoaoaoaoaoaoaoaoaEaoaYaaaAaGaqbGa}

  \item $\SSfamily{l+1}$ / $\CPSSfamily{k}{l+1}$ head iso sig:

    \texttt{jAMHwMMPkaabcddfgfghhiiTbgaTbdbaaRadbaaTbRaPbTbRafa}

  \item Move sequence:

    C,8 C,11 C,9 3,9 3,24 T,26 1,4 3,26 2,22 T,29 3,29 2,22 T,29 3,12 3,20 T,19 3,4 2,24 T,15 2,15 3,3 2,23 T,31 2,10 2,10 T,10 2,11 E,4 3,12 2,9 1,3 3,12 3,6 2,19 2,19 1,7 3,0 1,5 3,21 2,17 2,28 2,30 2,20 1,8 3,14 2,20 2,18 2,18 2,29 2,26 2,25 1,6 3,23 2,25 1,6 3,15 1,5 3,25 2,18 2,14 2,11 2,24 1,8 1,5 3,23 2,10 1,3 3,20 2,20 2,15 2,17 2,17 3,3 1,9 2,10 3,3 1,9 1,5 2,9 1,4

\end{itemize}

\section{Code}
\label{app:code}

\subsection*{4-Dimensional Connected Sum Script}

The script below includes Python code to generate the triangulation of $S^3 \times I$ called $\mathcal{C}$, as the function \texttt{double\_prism\_cylinder()}. 

\begin{lstlisting}[language=python]
from regina import *

# Construct a triangulation of S3xI where S3x{0} and S3x{1} are both 2-tetrahedron "pillow" 
#   triangulations of S3
# all gluings are performed such that only vertices with the same label are identified
def double_prism_cylinder():
  
  t = Triangulation4()
  for i in range(0,8):
    t.newPentachoron()

  # construct two disjoint tetrahedral prisms (pentachora 0-3 and pentachora 4-7)
  # end tetrahedra are face 4 of pentachora 0,4 and face 0 of pentachora 3,7
  for i in range(3):
    # note gluing is "diagonal" so that in pentachora i and i+4, faces 4-i and 3-i 
    # are always either glued in this step or will end up on the boundary
    t.pentachoron(i).join(3-i,t.pentachoron(i+1),Perm5())
    t.pentachoron(i+4).join(3-i,t.pentachoron(i+5),Perm5())

  # identify the two prisms along their "walls"
  for i in range(4):
    for j in range(0,3-i):
      t.pentachoron(i).join(j,t.pentachoron(i+4),Perm5())
    for j in range(5-i,5):
      t.pentachoron(i).join(j,t.pentachoron(i+4),Perm5())
  
  return t

# Puncture the triangulation by opening along the given face of the given pentachoron
# and inserting a double_prism_cylinder()
# pent_index is the integer index of a pentachoron, face_index is the index (0,1,2,3,4) of one of its faces
# this must not be a boundary face
# Modifies the existing triangulation
# Adds 8 new pentachora
# The new S3 boundary is a "pillow" given by face 0 of pentachora n+3 and n+7
def puncture_4d(t,pent_index,face_index):
  
  n = t.countFaces(4)
  t.insertTriangulation(double_prism_cylinder())
  
  # pop open along the given face of pent
  # glue in cylinder using the two tetrahedra of one of the S3 boundaries
  # face 4 (boundary) of pentachoron n is glued to pent
  # entry in the gluing table is the label of the given facet (i.e. ascending order)
  # face 4 (boundary) of pentachoron n+4 is glued to the adjacent pentachoron
  # entry in gluing table should be the same as the original entry for pent
  pent = t.pentachoron(pent_index)
  adj = pent.adjacentSimplex(face_index)
  gluing_to_pent = Perm5.rot((face_index-4) % 5)
  gluing_to_adj = pent.adjacentGluing(face_index) * Perm5.rot((face_index-4) % 5)
  
  pent.unjoin(face_index)
  t.pentachoron(n).join(4,pent,gluing_to_pent)
  t.pentachoron(n+4).join(4,adj,gluing_to_adj)

# Find the connected sum of two triangulations, by opening along the given faces of given pentachora
# in each triangulation, and inserting one end of a double_prism_cylinder() into each
# t1,t2: 4d triangulations
# pent1_index: int, index of a pentachoron in t1
# face1_index: ind, index (0,1,2,3,4) of one of its faces
# pent2_index: int, index of a pentachoron in t2
# face2_index: int, index (0,1,2,3,4) of one of its faces
# these faces must not be boundary faces
# sign: boolean, determines which of the two ways to perform the connected sum is used 
# NOTE: it is not guaranteed that the default sign=False will respect existing orientations
def connected_sum_4d(t1,t2,pent1_index,face1_index,pent2_index,face2_index,sign=False):
  
  n1 = t1.countFaces(4)
  n2 = t2.countFaces(4)
  t = Triangulation4()
  t.insertTriangulation(t1)
  t.insertTriangulation(t2)

  # puncture the copy of t1 along given face of pent1
  puncture_4d(t,pent1_index,face1_index)

  # pop open along the given face of pent2 in copy of t2
  # glue in punctured t1 using the two tetrahedra of the S3 boundary created by puncturing
  # similarly to in puncture_4d()
  # face 0 (boundary) of pentachoron n1+n2+3 is glued to pent2
  # entry in the gluing table is the label of the given facet (i.e. ascending order)
  # face 0 (boundary) of pentachoron n1+n2+7 is glued to the adjacent pentachoron
  # entry in gluing table should be the same as the original entry for pent2
  pent = t.pentachoron(n1+pent2_index)
  adj = pent.adjacentSimplex(face2_index)
  gluing_to_pent = Perm5.rot(face2_index % 5)
  gluing_to_adj = pent.adjacentGluing(face2_index) * Perm5.rot(face2_index % 5)
  
  # for the other connected sum, swap two vertices in each new gluing table entry
  if (sign):
    gluing_to_pent = gluing_to_pent * Perm5(0,2,1,3,4)
    gluing_to_adj = gluing_to_adj * Perm5(0,2,1,3,4)
  
  pent.unjoin(face2_index)
  t.pentachoron(n1+n2+3).join(0,pent,gluing_to_pent)
  t.pentachoron(n1+n2+7).join(0,adj,gluing_to_adj)

  return t
    # note gluing is "diagonal" so that in pentachora i and i+4, faces 4-i and 3-i 
    # are always either glued in this step or will end up on the boundary
    t.pentachoron(i).join(3-i,t.pentachoron(i+1),Perm5())
    t.pentachoron(i+4).join(3-i,t.pentachoron(i+5),Perm5())

  # identify the two prisms along their "walls"
  for i in range(4):
    for j in range(0,3-i):
      t.pentachoron(i).join(j,t.pentachoron(i+4),Perm5())
    for j in range(5-i,5):
      t.pentachoron(i).join(j,t.pentachoron(i+4),Perm5())
  
  return t

# Puncture the triangulation by opening along the given face of the given pentachoron
# and inserting a double_prism_cylinder()
# pent_index is the integer index of a pentachoron, face_index is the index (0,1,2,3,4) of one of its faces
# this must not be a boundary face
# Modifies the existing triangulation
# Adds 8 new pentachora
# The new S3 boundary is a "pillow" given by face 0 of pentachora n+3 and n+7
def puncture_4d(t,pent_index,face_index):
  
  n = t.countFaces(4)
  t.insertTriangulation(double_prism_cylinder())
  
  # pop open along the given face of pent
  # glue in cylinder using the two tetrahedra of one of the S3 boundaries
  # face 4 (boundary) of pentachoron n is glued to pent
  # entry in the gluing table is the label of the given facet (i.e. ascending order)
  # face 4 (boundary) of pentachoron n+4 is glued to the adjacent pentachoron
  # entry in gluing table should be the same as the original entry for pent
  pent = t.pentachoron(pent_index)
  adj = pent.adjacentSimplex(face_index)
  gluing_to_pent = Perm5.rot((face_index-4) % 5)
  gluing_to_adj = pent.adjacentGluing(face_index) * Perm5.rot((face_index-4) % 5)
  
  pent.unjoin(face_index)
  t.pentachoron(n).join(4,pent,gluing_to_pent)
  t.pentachoron(n+4).join(4,adj,gluing_to_adj)

# Find the connected sum of two triangulations, by opening along the given faces of given pentachora
# in each triangulation, and inserting one end of a double_prism_cylinder() into each
# t1,t2: 4d triangulations
# pent1_index: int, index of a pentachoron in t1
# face1_index: ind, index (0,1,2,3,4) of one of its faces
# pent2_index: int, index of a pentachoron in t2
# face2_index: int, index (0,1,2,3,4) of one of its faces
# these faces must not be boundary faces
# sign: boolean, determines which of the two ways to perform the connected sum is used 
# NOTE: it is not guaranteed that the default sign=False will respect existing orientations
def connected_sum_4d(t1,t2,pent1_index,face1_index,pent2_index,face2_index,sign=False):
  
  n1 = t1.countFaces(4)
  n2 = t2.countFaces(4)
  t = Triangulation4()
  t.insertTriangulation(t1)
  t.insertTriangulation(t2)

  # puncture the copy of t1 along given face of pent1
  puncture_4d(t,pent1_index,face1_index)

  # pop open along the given face of pent2 in copy of t2
  # glue in punctured t1 using the two tetrahedra of the S3 boundary created by puncturing
  # similarly to in puncture_4d()
  # face 0 (boundary) of pentachoron n1+n2+3 is glued to pent2
  # entry in the gluing table is the label of the given facet (i.e. ascending order)
  # face 0 (boundary) of pentachoron n1+n2+7 is glued to the adjacent pentachoron
  # entry in gluing table should be the same as the original entry for pent2
  pent = t.pentachoron(n1+pent2_index)
  adj = pent.adjacentSimplex(face2_index)
  gluing_to_pent = Perm5.rot(face2_index % 5)
  gluing_to_adj = pent.adjacentGluing(face2_index) * Perm5.rot(face2_index % 5)
  
  # for the other connected sum, swap two vertices in each new gluing table entry
  if (sign):
    gluing_to_pent = gluing_to_pent * Perm5(0,2,1,3,4)
    gluing_to_adj = gluing_to_adj * Perm5(0,2,1,3,4)
  
  pent.unjoin(face2_index)
  t.pentachoron(n1+n2+3).join(0,pent,gluing_to_pent)
  t.pentachoron(n1+n2+7).join(0,adj,gluing_to_adj)

  return t
\end{lstlisting}

\subsection*{Code to Generate the Families}

\begin{lstlisting}[language=python]
from regina import *

# P(k) = #^k CP^2
def cp2k(k):
  if (k<=0):
    return Triangulation4("cAkaabb+aoa+aoa")
  elif k==1:
    return Triangulation4("eAMMcaabccdd+aoa+aAaqbqbGa")

  first = Triangulation4("eAMMcaabccdd+aoa+aAaqbqbGa")
  first.removePentachoron(first.pentachoron(3))
  unit_cp2 = Triangulation4("eAMMcaabccdd+aoa+aAaqbqbGa")
  unit_cp2.removePentachoron(unit_cp2.pentachoron(3))
  unit_cp2.removePentachoron(unit_cp2.pentachoron(0))
  last = Triangulation4("eAMMcaabccdd+aoa+aAaqbqbGa")
  last.removePentachoron(last.pentachoron(0))

  t = Triangulation4()
  t.insertTriangulation(first)
  for i in range(k-2):
    t.insertTriangulation(unit_cp2)
    t.pentachoron(t.countPentachora()-3).join(4,t.pentachoron(t.countPentachora()-2),Perm5(1,2,0,3,4))
  t.insertTriangulation(last)
  t.pentachoron(t.countPentachora()-4).join(4,t.pentachoron(t.countPentachora()-3),Perm5(1,2,0,3,4))

  return t

# A(k) = #^k (CP^2 # -CP^2)
def cp2PMk(k):
  if (k<=0):
    return Triangulation4("cAkaabb+aoa+aoa")
  elif k==1:
    return Triangulation4("gALAMQaacdcdeeff+aoawbaaYa+ayaYa+aoa")

  first = Triangulation4("gALAMQaacdcdeeff+aoawbaaYa+ayaYa+aoa")
  first.removePentachoron(first.pentachoron(5))
  unit_cp2PM = Triangulation4("gALAMQaacdcdeeff+aoawbaaYa+ayaYa+aoa")
  unit_cp2PM.removePentachoron(unit_cp2PM.pentachoron(5))
  unit_cp2PM.removePentachoron(unit_cp2PM.pentachoron(0))
  last = Triangulation4("gALAMQaacdcdeeff+aoawbaaYa+ayaYa+aoa")
  last.removePentachoron(last.pentachoron(0))

  t = Triangulation4()
  t.insertTriangulation(first)
  for i in range(k-2):
    t.insertTriangulation(unit_cp2PM)
    t.pentachoron(t.countPentachora()-6).join(4,t.pentachoron(t.countPentachora()-4),Perm5())
  t.insertTriangulation(last)
  t.pentachoron(t.countPentachora()-7).join(4,t.pentachoron(t.countPentachora()-5),Perm5())

  return t

# E(k) = #^k S^2 x S^2
def s2xs2k(k):
  if (k<=0):
    return Triangulation4("cAkaabb+aoa+aoa")
  elif k==1:
    return Triangulation4("gALAMQaacdcdeeff+aoawbaaGa+ayaGa+aoa")

  first = Triangulation4("gALAMQaacdcdeeff+aoawbaaGa+ayaGa+aoa")
  first.removePentachoron(first.pentachoron(5))
  unit_s2xs2 = Triangulation4("gALAMQaacdcdeeff+aoawbaaGa+ayaGa+aoa")
  unit_s2xs2.removePentachoron(unit_s2xs2.pentachoron(5))
  unit_s2xs2.removePentachoron(unit_s2xs2.pentachoron(0))
  last = Triangulation4("gALAMQaacdcdeeff+aoawbaaGa+ayaGa+aoa")
  last.removePentachoron(last.pentachoron(0))

  t = Triangulation4()
  t.insertTriangulation(first)
  for i in range(k-2):
    t.insertTriangulation(unit_s2xs2)
    t.pentachoron(t.countPentachora()-6).join(4,t.pentachoron(t.countPentachora()-4),Perm5())
  t.insertTriangulation(last)
  t.pentachoron(t.countPentachora()-7).join(4,t.pentachoron(t.countPentachora()-5),Perm5())

  return t

# D(k,l) = (#^k CP^2) # (#^l S^2 x S^2)
def cp2k_s2xs2l(k,m):
  if m<=0:
    return cp2k(k)
  elif k<=0:
    return s2xs2k(m)

  first = Triangulation4("eAMMcaabccdd+aoa+aAaqbqbGa")
  first.removePentachoron(first.pentachoron(3))
  unit_cp2 = Triangulation4("eAMMcaabccdd+aoa+aAaqbqbGa")
  unit_cp2.removePentachoron(unit_cp2.pentachoron(3))
  unit_cp2.removePentachoron(unit_cp2.pentachoron(0))
  unit_s2xs2 = Triangulation4("gALAMQaacdcdeeffPbgaVbaaJafbyaJafbva")
  unit_s2xs2.removePentachoron(unit_s2xs2.pentachoron(5))
  unit_s2xs2.removePentachoron(unit_s2xs2.pentachoron(0))
  last = Triangulation4("gALAMQaacdcdeeffPbgaVbaaJafbyaJafbva")
  last.removePentachoron(last.pentachoron(0))

  t = Triangulation4()
  t.insertTriangulation(first)
  for i in range(k-1):
    t.insertTriangulation(unit_cp2)
    t.pentachoron(t.countPentachora()-3).join(4,t.pentachoron(t.countPentachora()-2),Perm5(1,2,0,3,4))
  if m==1:
    t.insertTriangulation(last)
    t.pentachoron(t.countPentachora()-6).join(4,t.pentachoron(t.countPentachora()-5),Perm5(1,2,0,4,3))
  else:
    t.insertTriangulation(unit_s2xs2)
    t.pentachoron(t.countPentachora()-5).join(4,t.pentachoron(t.countPentachora()-4),Perm5(1,2,0,4,3))
    for i in range(m-2):
      t.insertTriangulation(unit_s2xs2)
      t.pentachoron(t.countPentachora()-6).join(3,t.pentachoron(t.countPentachora()-4),Perm5())
    t.insertTriangulation(last)
    t.pentachoron(t.countPentachora()-7).join(3,t.pentachoron(t.countPentachora()-5),Perm5())

  return t
\end{lstlisting}

\end{appendices}

\end{document}